%
%

\magnification=1200

\font\titfont=cmr10 at 12 pt


\font\AAA=cmr14 at 12pt
\font\BBB=cmr14 at 8pt

\overfullrule=0in

\def\smfrac#1#2{\hbox{${#1\over #2}$}}
\def\oa#1{\overrightarrow #1}
\def\dim{{\rm dim}}

\def\deg{{\rm deg}}

\def\Hess{{\rm Hess}}
\def\Sym{\rm Sym}

\def\tr{{\rm tr}}

\def\span{{\rm span\,}}

\def\End{{\rm End}}

\def\arr{\longrightarrow}
\def\supp{{\rm supp}}

\def\foral{\qquad {\rm for\ all\ \ }}
\def\fpsh{{\rm PSH}(X,\f)}
\def\Core{{\rm Core}}


\def\Theorem#1{\medskip\noindent {\AAA T\BBB HEOREM \rm #1.}}
\def\Prop#1{\medskip\noindent {\AAA P\BBB ROPOSITION \rm  #1.}}
\def\Cor#1{\medskip\noindent {\AAA C\BBB OROLLARY \rm #1.}}
\def\Lemma#1{\medskip\noindent {\AAA L\BBB EMMA \rm  #1.}}
\def\Remark#1{\medskip\noindent {\AAA R\BBB EMARK \rm  #1.}}
\def\Note#1{\medskip\noindent {\AAA N\BBB OTE \rm  #1.}}
\def\Def#1{\medskip\noindent {\AAA D\BBB EFINITION \rm  #1.}}

\def\Ex#1{\medskip\noindent {\AAA E\BBB XAMPLE \rm    #1.}}

\def\pf{\medskip\noindent {\bf Proof.}\ }
\def\qed{\hfill  $\vrule width5pt height5pt depth0pt$}

\def\qedqed{\hfill  $\vrule width5pt height5pt depth0pt$ $\vrule width5pt height5pt depth0pt$}

\def\hk{\_{\rm l}\,}

   \def\cp{{\cal P}}
   
\def\ce{{\cal E}}   
\def\ch{{\cal H}}   \def\cm{{\cal M}}
   
\def\cd{{\cal D}}

\def\cp{{\cal P}}

\def\ccr{{\cal  R}}

\def\wt{\widetilde}
\def\wh{\widehat}

\def\and{\qquad {\rm and} \qquad}
\def\arr{\longrightarrow}

\def\bbr{{\bf R}}\def\bbh{{\bf H}}

\def\bbz{{\bf Z}}

\def\a{\alpha}
\def\b{\beta}
\def\d{\delta}
\def\e{\epsilon}
\def\f{\phi}

\def\o{\omega}

\def\s{\sigma}
\def\x{\xi}

\def\D{\Delta}
\def\L{\Lambda}
\def\G{\Gamma}
\def\O{\Omega}

\def\bd{\partial}
\def\bdf{\partial_{\f}}
\def\fp{$\phi$-plurisubharmonic }
\def\fh{$\phi$-pluriharmonic }

\def\ffl{$\f$-flat}
\def\PH#1{\widehat {#1}}
\def\psh{plurisubharmonic \ }

\def\lp{\Lambda_+(\f)}
\def\lpp{\Lambda^+(\f)}
\def\bo{\partial \Omega}
\def\fc{$\phi$-convex }
\def\PSH{\rm PSH}

\def\BM{\lambda}

\def\RH{\overline{\ch}^\f }

\def\YY{1}
\def\AA{1}

\def\DD{2}
\def\EE{3}

\def\GG{6}
\def\HH{4}
\def\II{5}

\ 
\vskip .3in

\centerline{\titfont DUALITY OF POSITIVE CURRENTS}\smallskip
\centerline{\titfont  AND PLURISUBHARMONIC 
FUNCTIONS  }
\smallskip

\centerline{\titfont IN CALIBRATED GEOMETRY }
\bigskip

\centerline{\titfont F. Reese Harvey and H. Blaine Lawson, Jr.$^*$}
\vglue .9cm
\smallbreak\footnote{}{ $ {} \sp{ *}{\rm Partially}$  supported by
the N.S.F. }

\vskip .3in
\centerline{\bf ABSTRACT } \medskip
  \font\abstractfont=cmr10 at 10 pt

{{\parindent= .6in\narrower\abstractfont \noindent
Recently the authors showed that there is a robust potential
theory attached to any calibrated manifold $(X,\f)$. In particular,
 on $X$ there exist  $\f$-plurisubharmonic functions, $\f$-convex 
domains, $\f$-convex boundaries, etc.,  all inter-related 
and having a  number of good properties.
In this paper we show that, in a strong sense, the plurisubharmonic
functions are the polar duals of the $\f$-submanifolds, or more generally,
the   $\f$-currents studied in the original paper on calibrations.
In particular, we establish an analogue of Duval-Sibony Duality
which characterizes points in the $\f$-convex hull of a compact set
$K\subset X$ in terms of $\f$-positive Green's currents on $X$  and Jensen 
measures on $K$.  We also characterize boundaries of $\f$-currents
entirely in terms of $\f$-plurisubharmonic functions.
Specific calibrations are used as examples throughout.
Analogues of the Hodge Conjecture in  calibrated geometry are considered.

}}

\vfill\eject\
\vskip 1in

\centerline{\bf TABLE OF CONTENTS} \bigskip

\qquad \YY. Introduction.\smallskip

\qquad \DD.  Positive Currents in Calibrated Geometries. \smallskip

\qquad \EE.    Duality with Plurisubharmonic Functions.\smallskip

\qquad \HH. Boundary Duality.\smallskip

\qquad \II. $\f$-Flat Hypersurfaces and Functions which are $\f$-Pluriharmonic mod $d$.\smallskip

\qquad \GG.  Hodge Manifolds.\smallskip

\qquad\qquad\qquad\qquad Appendix:  {The Reduced $\f$-Hessian. }  \smallskip

\vfill\eject

\centerline{\bf \YY. Introduction.}
\medskip

Calibrated geometries are  geometries
of distinguished minimal varieties  determined by a fixed, closed differential form $\f$ on
a riemannian manifold $X$.
A basic example is that of a K\"ahler manifold (or more generally a symplectic manifold,
with compatible complex structure) where the distinguished submanifolds are the  holomorphic curves.
However, there exist many other interesting  geometries, each carrying a wealth of $\f$-submanifolds, particularly  on spaces with special holonomy.  These have attracted particular attention in recent
years due to their appearance in generalized Donaldson theories and in modern versions of string theory in Physics.

Recently   it was shown  [HL$_5$] that a surprisingly large part of classical
pluripotential theory can be carried over to any calibrated manifold $(X,\f)$. In particular,
there is a natural family $\PSH(X,\f)$ of $\f$-\psh functions possessing essentially
all of the important properties of the plurisubharmonic functions in complex analysis.
For example, the family is closed under composition with convex increasing functions,
under taking the maximum of elements in the family, and under taking upper envelopes.
Furthermore, there is a notion of  $\f$-plurisubharmonic distribution, and if $\f$ is {\sl elliptic} 
(i.e., $\f$ involves all the variables at each point -- see below),
then the cone of these  is weakly closed, and every such distribution has a unique upper semi-continuous, $L^1_{\rm loc}$-representative.

There is   a  notion  of $\f$-convexity for   $(X,\f)$, defined in terms 
of  the \fp hulls of  compact subsets, and  characterized
by the existence of a proper \fp exhaustion function.   
There is also a notion of boundary convexity for  
 domains
$\O\subset\subset X$, and various aspects of the Levi problem hold.
For example, if $\bo$ is strictly $\f$-convex, then so is $\O$ 
(provided $X$ admits some strictly psh-function).
Furthermore, on strictly convex domains, the Perron-Bremermann method can be applied to solve the
Dirichlet problem for  $\f$-partially pluriharmonic functions --- the analogue of the
homogeneous Monge-Amp\`ere equation on $X$,  [HL$_6$]. These methods also lead
to various notions of capacity in calibrated geometry.

Now it is known that   \fp functions and   $\f$-submanifolds are intimately related.  
For example:   
\Theorem{\YY.1.  ([HL$_5$])}
{\sl The restriction of a \fp function to a $\f$-submanifold $M$ is subharmonic in the induced riemannian metric on $M$.}

\smallskip 

The main thrust of this paper could be summarized in the following 
\smallskip\noindent
{\bf General Principle:}  {\sl The \fp functions and $\f$-submanifolds (more specifically, the positive $\f$-currents) are polar duals of one-another.}
\smallskip

We shall see this in \S \EE\ where points $x$ in the $\f$-convex hull of a compact set $K\subset X$ are characterized by the existence of a $\f$-positive {\sl Green's current  }$T$ and a Jensen measure
$\mu$ on $K$ satisfying the generalized {\sl Poisson-Jensen equation}
$$
 \partial_\f \partial T \ =\ \mu - [x].
$$
This is an extension of Duval-Sibony duality [DS]  to general calibrated geometries.

We shall also see this in \S \HH \  where the boundaries of $\f$-positive currents are characterized.
 \medskip

 For the reader's sake, here in the introduction, we  now briefly review the basic definitions and results from
 [HL$_5$] to which this paper is a sequel. A {\sl calibrated manifold} is a pair
 $(X,\f)$ where $X$ is a riemannian manifold and $\f$ is an exterior $p$-form
 on $X$ which is $d$-closed and has comass one, i.e.,
 $$
 \f(\x)  \leq \ 1
 $$
for every unit simple $p$-vector $\x$ on $X$. The 
{\sl $\f$-Grassmannian} is the set of {\sl $\f$-planes:}
$$
G(\f)\ \equiv\ \{\x   :  \f(\x)=1\}.
$$

The \fp functions on  $(X,\f)$ are defined by a  second order
differential operator
$
\ch^\f : C^\infty(X) \ \to\ \ce^p(X),
$
called the {\sl $\f$-Hessian}, defined by  
$$
\ch^\f(f)\ =\ \BM_\f(\Hess f)
$$
where $\Hess f$ is the riemannian hessian of $f$ and $\BM_\f:\End(TX) \to \L^pT^*X$
is the bundle map given by  $\BM_\f(A) = D_{A^*}(\f)$ where  $D_{A^*}:\L^pT^*X\arr \L^pT^*X$ is the natural extension of $A^*:T^*X \to T^*X$ as a derivation.
When $\nabla\f=0$,   there is a natural factorization
$$
\ch^\f \ =\ d d^\f
$$
where  $d^\f: C^\infty(X) \ \to\ \ce^{p-1}(X)$ is given by
$$
d^\f f\ \equiv \ \nabla f \hk \f.
$$
In general these operators are related by the equation: $\ch^\f  f = dd^\f f - \nabla_{\nabla f} (\f)$.

A  function
$f\in C^\infty(X)$ is defined to be {\bf  \fp} if 
$$
\ch^\f(f)(\x) \geq 0\qquad {\rm for\  all\ }  \x\in G(\f)
$$
The function $f$ is called  {\sl strictly \fp } at a point $x\in X$ if $\ch^\f(f)(\x) > 0$
 for all $\f$-planes $\x$ at $x$.  
 It is $\f$-{\sl pluriharmonic} if $\ch^\f(f)(\x)=0$ for all $\phi$-planes $\x$.
We denote by $\PSH(X,\f)$ the convex cone of smooth \fp functions on $X$.

When $X$ is a complex manifold with a K\"ahler form $\o$, one easily computes that
$d^\o = d^c$, the conjugate differential.  In this case, $\ch^\o = dd^\o =dd^c$ and the 
$\o$-planes correspond to the complex lines in $TX$. Hence, the definitions above coincide
with the classical notions of plurisubharmonic and pluriharmonic functions on $X$.
 
A fundamental property of the $\f$-Hessian is that for any $\f$-plane $\x$, one has
$$
\left(\ch^\f f\right)(\x)\ =\ {\rm trace}\left\{ \Hess f\bigr|_{\x}  \right\}.
$$

The first concept addressed in  [HL$_5$]   is the analogue
of pseudoconvexity in complex geometry.
 Let  $(X,\f)$ be a calibrated manifold and $K\subset X$ a closed
subset. By the {\bf $\f$-convex hull of $K$ } we mean the subset
$$
\wh K\ =\ \{x\in X : f(x) \leq \sup_K f \ \ {\rm for \ all \  } f\in \PSH(X,\f)\}
$$
The manifold $(X,\f)$ is said to be {\bf  $\f$-convex} if $K\subset \subset X \ \Rightarrow
\  \wh K \subset \subset X$  for all $K$.

\Theorem {\YY.2.  ([HL$_5$])} {\sl A calibrated manifold  $(X,\f)$ is $\f$-convex if and only if 
it admits a \fp proper exhaustion function $f:X\to \bbr$.} \medskip

The manifold $(X,\f)$ will be called {\sl strictly $\f$-convex} if it admits an exhaustion function $f$ which is  strictly $\f$-plurisubharmonic, 
and it will be called  {\sl strictly $\f$-convex at infinity}  if $f$ is strictly \fp  outside of a compact subset.

Note that in complex geometry, \fc manifolds are Stein, and manifolds which are \fc at infinity are 
 strongly pseudoconvex. In the latter case there is a distinguished compact subset
 consisting of exceptional subvarieties. This set has the following general analogue.

The {\bf core} of $X$   is   the set of points
$x\in X$ with the property that no $f\in\PSH(X,\f)$ is strictly \fp at $x$. \smallskip


\Theorem {\YY.3. ([HL$_5$])} {\sl   Suppose $X$ is $\f$-convex.  Then  $X$ is strictly $\f$-convex at infinity if and only if {\rm Core}$(X)$ is compact, and $X$ is strictly $\f$-convex if and only if {\rm Core}$(X)=\emptyset$.}
\smallskip

A very general construction of strictly \fc manifolds is given  in  [HL$_5$, \S \GG].
   The construction is based on the notion of {\sl $\phi$-free submanifolds}, which directly
              generalize the totally real submanifolds in complex geometry.

We now sketch the contents of this paper.
\medskip
\centerline
{$\f$\AAA -P\BBB OSITIVE \AAA C\BBB URRENTS.}\smallskip
 
In \S \DD\ we review the theory of $\f$-positive currents originally introduced in [HL$_3$].
We recall that a $p$-dimensional current $T$ is called {\bf $\f$-positive }
if it is representable by integration and its generalized tangent $p$-vector
$$
\overrightarrow{T} \in {\rm Convex Hull }(G(\f))\qquad \|T\|-a.e.
$$
where $\|T\|$ denotes the total variation measure of $T$. Examples include $\f$-submanifolds
and, more generally, rectifiable $\f$-currents. By Almgren's   Theorem [A] we know that 
rectifiable $\f$-currents $T$ with $dT=0$ are regular, that is, given by integration over 
$\f$-submanifolds with positive integer multiplicities, outside a closed subset of Hausdorff
dimension $p-2$.

$\f$-Positive currents generalize the positive currents in complex geometry, and $d$-closed
rectifiable $\f$-currents generalize positive holomorphic chains.

If $T$ is a $\f$-positive current (with compact support), then
$$
\supp\, T\ \subset \ \wh{\supp ( d T)} \cup \Core(X).
\eqno{(\YY.1)}
$$
In particular, if $dT=0$, then $\supp T \subset   \Core(X)$, and if $X$ is strictly $\f$-convex 
(i.e., $\Core(X)=\emptyset$), then there exist no $d$-closed $\f$-positive currents
 with compact support on $X$.

 Section \DD \   summarizes the known facts concerning  $\f$-positive currents.
These include compactness theorems, regularity theorems, mass-minimizing properties, monotonicity
properties, 
and dual characterizations.

\medskip
\centerline
{\AAA T\BBB HE \AAA  	S\BBB  UPPORT  \AAA  	L\BBB  EMMA.}\smallskip
\smallskip

Assume for now that the calibration  $\f$ is parallel, and consider the adjoint of the operator $dd^\f:\ce^0(X)\to \ce^p(X)$ which can be written as 
$$
\partial_\f\partial:\ce'_p(X) \ \arr\ \ce'_0(X) 
$$
where $\partial:\ce'_p(X) \to \ce'_{p-1}(X) $ denotes the usual adjoint of $d:\ce^{p-1}(X)\to \ce^p(X)$ and 
$\partial_\f:\ce'_{p-1}(X) \ \arr\ \ce'_0(X) $ is the adjoint of $d^\f$, 
defined by $(\partial_\f R)(f)\equiv R( d^\f f)$.

Positive currents $T$ with the property: $\partial_\f\partial T\leq 0$ (i.e.,  $\partial_\f\partial T$ is a non-positive measure), satisfy a version of (\YY.1) above.

\Lemma{\EE.2} {\sl Suppose $T$ is a $\f$-positive current with compact support on $X$ which satisfies 
$\partial_\f\partial T\leq 0$ outside a compact subset $K\subset X$. Then
$$
\supp\, T\ \subset \ \wh K \cup \Core(X).
$$
In particular, if $\partial_\f\partial T\leq 0$ on $X$, then  $\supp \,T \subset \Core(X)$.}\medskip

Another consequence is the following.  Suppose $X$ is strictly  $\f$-convex. If
$T$ is  a  $\f$-positive current with 
 $ \partial_\f\partial T \leq 0$ on $X-K$, then $\supp\, T\subset \wh K$.
In fact, it turns out that  the  points   $ x\in \wh K$ can be characterized in terms
of certain $\f$-positive currents $T$ which satisfy $\partial_\f\partial T = -[x]$ in  $X-K$.
We  discuss this next in greater detail.

\medskip
\centerline
{\AAA D\BBB UVAL - \AAA S\BBB IBONY \AAA D\BBB UALITY.} 
\smallskip

Points in the $\f$-convex hull of a compact set $K\subset X$ have a useful characterization
in terms of $\f$-positive currents and certain Poisson-Jensen measures.  The following
 results generalize work of Duval and Sibony [DS] in the complex case. 
 They remain valid (as does Lemma \EE.2 above) when $\f$ is not parallel if  the operator $\partial_\f\partial$ is replaced with $\ch_\f$.

Let $K\subset X$ be a compact subset and $x$ a point in $X-K$.  By a {\bf Green's current
for } $(K,x)$ we mean a $\f$-positive current $T$ which satisfies
$$
\partial_\f\partial T = \mu-[x]
\eqno{(\YY.2)}
$$
where $\mu$ is a probability measure with support on $K$ and $[x]$ denotes the Dirac measure
at $x$. In this case $\mu$ is called a {\bf Poisson-Jensen} measure for $(K,x)$.
By the remarks above we see that $x\in \wh K$.  In fact we have the following.

\Theorem{\EE.8}  {\sl  Suppose $\f$ is parallel and $X$ is strictly $\f$-convex.
Let $K\subset X$ be a compact subset and $x\in X-K$.  Then there exists a Green's current
for $(K,x)$ if and only if $x\in \wh K$.}

\medskip
We note that if $M\subset X$ is a compact $\f$-submanifold with boundary, and if $G_x$ 
is the Green's function for the riemannian laplacian on $M$ with singularity at $x\in M-\partial M$,
then $\partial_\f\partial(G_x [M]) = \mu-[x]$ for a probability measure $\mu$ on $\partial M$.

\medskip
As a application we obtain the following approximation result.
A domain $\O\subset X$ is said to be {\bf \fc relative to $X$} if 
$K\subset\subset \O \ \Rightarrow  \wh{K}_X\subset\subset \O$.
\Prop{\EE.16} {\sl  Suppose $\f$ is parallel and $X$ is strictly $\f$-convex.
An open subset $\O\subset X$ is  \fc relative to $X$ if and only if 
$\PSH(X,\f)$ is dense in $\PSH(\O,\f)$.}

 \vfill\eject
\centerline
{\AAA B\BBB OUNDARY \AAA D\BBB UALITY.}
\smallskip

A very natural question in calibrated geometry is the following:  Given a compact
oriented submanifold $\G\subset X$ of dimension $p-1$, when does there exist
a compact $\f$-submanifold $M$ with $ \partial M=\G$? A companion question is: Given 
a compactly supported current $S$ of dimension $p-1$ in $X$, when does there exist
a compactly supported $\f$-positive current $T$ with $\partial T=\G$? 
For this second question there is a complete answer
when   $X$ is strictly \fc and $\f$ is exact.

\Theorem {\HH.1} {\sl   Fix $(X,\f)$ as above, and consider  a current
$S\in \ce'_{p-1}(X)$. Then $S=\partial T$  for some  $\f$-positive current $T\in \ce'_p(X)$
if and only if }
$$
\int_S \, \a\ \geq \ 0\qquad\ \ {\sl for\ all\ } \a\in \ce^{p-1}(X) \ \ \ {\sl such\ that\ } d\a \ {\sl is\ }
\lpp{\sl -positive}
$$
\Note{} $\lpp$-positive means that $d\a(\x) \geq 0$ for all $\x\in G(\f)$.
\smallskip

There is a similar result for compact   manifolds $(X,\f)$ with
no condition on $\f$.  
\def\C{\lambda}

\Theorem{\HH.3} {\sl Suppose $(X, \f)$ is a compact calibrated manifold.  Fix  $S\in \ce'_{p-1}(X)$
and $\C>0$.    Then  the following are equivalent.  
\smallskip
(i) \ There  exists a $\lp$-positive current $T\in \ce'_p(X)$ with $S=\partial T$  and    
$M(T)\leq \C$.  \smallskip

(ii)\  $\int_S\a \ \geq\ -\C$ for all $\a \in \ce^{p-1}(X)$ such that $d\a+\f$ is  $\lpp$-positive.}

\medskip
There is also a result for non-compact boundaries on general manifolds $(X,\f)$ (Theorem \HH.4 below).
Versions of these theorems in K\"ahler geometry appear in [HL$_4$].

 \medskip
\centerline
{\AAA $\f$-F\BBB LAT  \AAA H\BBB YPERSURFACES.} 
\smallskip

In Section \II \ we expand the notion of $\f$-pluriharmonic functions to include functions $f$
which are $\f$-pluriharmonic modulo the ideal generated by $df$.
In most interesting geometries these functions are characterized by the fact
that their level sets are {\bf $\f$-flat}, i.e., the trace of the second fundamental form
on all tangential $\f$-planes is zero.  These functions are important for the boundary problem.
If $f$ is such a function defined in a neighborhood of a compact $\f$-submanifold with boundary  
 $M\subset X$, then
$$
\inf_{\partial M} f \ \leq \  f(x) \ \leq \ \sup_{\partial M} f \qquad\ \ {\rm for \ \ } x\in M.
$$

\medskip
\centerline
{\AAA G\BBB ENERALIZED \AAA H\BBB ODGE  \AAA M\BBB ANIFOLDS.}\smallskip
\smallskip

 In Section \GG\ we discuss analogues of  Hodge manifolds in the general calibrated setting.
We also examine various analogues of the Hodge Conjecture in these spaces.

\medskip      
{\sl Note.}  For simplicity we always assume the manifold $X$ to be oriented.
Following de Rham we denote by  $\ce^p(X)$ the 
space of smooth differential forms of degree $p$ on $X$,  and by  $\ce_p'(X)$
its topological dual space of currents of dimension $p$ with compact support 
on $X$
\medskip

 The authors would like to thank Robert Bryant for useful comments and conversations
related to this paper.

\vfill\eject


\vskip .3in

\centerline{\bf \DD. Positive Currents in Calibrated Geometries.}

\medskip

The important classical notion of  a positive current on a complex manifold has an analogue
on any calibrated manifold.  This concept was introduced   in section II
of [HL$_3$].  In this Section we review that  material with some of the terminology and notation updated.  

Suppose $\f$ is a calibration on a riemannian manifold $X$.  Let $p=\deg \f$ and $n=\dim X$.
The $\f$-Grassmannian, denoted $G(\f)$, consists of the unit simple vectors $\x \in G(p,TX)\subset \L_pTX$ with $\f(\x)=1$, i.e.,  the $\f$-planes.

On a calibrated manifold $(X,\f)$ the three concepts we wish to discuss are:
\smallskip

\qquad  a)\ \ $\f$-submanifolds,\smallskip

\qquad  b)\ \ rectifiable $\f$-currents, and $\f$-cycles\smallskip

\qquad   c)\ \ $\f$-positive  currents.\smallskip

A $\f$-submanifold is, of course, a smooth oriented submanifold $M$ whose oriented
tangent space is a $\f$-plane at every point, i.e., $\oa M_x \equiv \oa T_x M\in G(\f)$ for all $x\in M$.

Suppose $T$ is a locally rectifiable $p$-dimensional current (cf. [F$_1$]) on $X$. Then its generalized tangent space  is a unit simple vector $\oa T\in G(p, TX)$  
at $\|T\|$ almost every point, where $\|T\|$ denotes the generalized 
volume measure associated with $T$.

\Def{\DD.1}  A {\bf  rectifiable $\f$-current}  is a locally rectifiable current
$T$ with $\oa T \in G(\f)$ for $\|T\|$- a.a. points in $X$.  A {\bf  $\f$-cycle}
is a rectifiable $\f$-current which is $d$-closed.

\medskip

\Remark{\DD.2} We shall see below (Theorem \DD.10) that $\f$-cycles always have a particularly nice
local  structure.  The  strongest  result of this kind occurs in the K\"ahler case, with $\f=\o^p/p!$,
where  a theorem of  J. King [K] states that each 
$\f$-cycle is a positive holomorphic cycle, i.e., a locally finite sum of $p$-dimensional complex analytic subvarieties  with positive integer coefficients.  
Results about the singular structure in the special Lagrangian case have been obtained by
  Joyce, Haskins, Kapouleas, Pacini, and others (cf.  [J$_*$], [Ha$_*$], [HaK], [HaP]).

\Remark{\DD.3}  On a general calibrated
manifold $(X,\f)$ one can also consider $d$-closed rectifiable currents $T$ with $\pm \overrightarrow T\in G(\f)$ for $\|T\|$-a.a. points.  In the K\"ahler case $T$ must be a holomorphic chain by results in [HS], [S] and [Alex].  However, nothing is known about the structure of such currents for any of the other standard calibrations.  

\medskip

An understanding of the definition of a $\f$-positive current is a little more complicated.

Recall (cf.  [F$_1$]) that a current $T$ is {\bf representable by integration} if $T$ has measure coefficients when expressed as a generalized differential form.  Equivalently, the mass norm $M_K(T)$ of $T$ on each compact set $K$, is finite.  Associated with such a current $T$ is a Radon measure $\|T\|$  and a generalized tangent space $\oa T_x\in\L_pT_xX$ defined for $\|T\|$-a.a. points $x$. Recall that each $\oa T_x$ has mass norm one.  For any $p$-form $\a$ with compact support
$$
T(\a)\ =\ \int\a(\oa T)\, d\|T\|
\eqno{(\DD.1)}
$$

\Def{\DD.4}  At each point $x\in X$ let $\L(\f)$ denote the span of $G(\f)\subset \L_p TX$, 
and let 
$$
\lp \ \subset \ \L(\f)
$$
denote the convex cone  on $G(\f)$ with vertex the origin.  The $p$-vectors
$\x\in \lp$ will be called {\bf $\lp$-positive}.
\medskip

Note that $\lp$ is just the cone on $ {\rm ch}  \, G(\f)$, the convex hull of the $\f$-Grassmannian.

The following Lemma is needed for a robust understanding of the definition of a 
$\f$-positive current.

\Lemma{\DD.5} {\sl The following conditions are equivalent:
\smallskip

\qquad\qquad\qquad 1)\ \ $\oa T \in \L_+(\f) \qquad\ \  \|T\| $-a.e.
\smallskip

\qquad\qquad\qquad 2)\ \ $\oa T \in {\rm ch}  \, G(\f) \qquad \|T\|$-a.e.
\smallskip

\qquad\qquad\qquad 3)\ \ $\phi(\oa T )=1 \qquad\ \ \ \   \|T\|$-a.e.
\medskip\noindent
}

The proof is provided later.

\Def{\DD.6}  A {\bf $\f$-positive  current}  is a $p$-dimensional current $T$
which is representable by integration and for which  the equivalent 
conditions of Lemma \DD.5 are satisfied.\medskip

By Condition 3) of \DD.5  and (\DD.1)
$$
  T(\f)\ =\ \int \f(\oa T)\, d\|T\| \ = M(T)
\eqno{(\DD.2)}
$$
for all $\f$-positive currents $T$.  This fact  has important implications.
\medskip

\Prop{\DD.7} {\sl  Suppose $T$ is a compactly supported $p$-dimensional current which is representable by integration.  Then
$$
T(\f)\ \leq \ M(T)
$$
with equality if and only if $T$ is a $\f$-positive current.

Consequently, any $\f$-positive current $T_0$ with compact support is homologically mass-minimizing, i.e., 
$$
M(T_0)\ \leq \ M(T)
\eqno{(\DD.3)}
$$
for any $T=T_0+dS$ where $S$ is a $(p+1)$-dimensional current with compact support.  Furthermore,
equality holds in (\DD.3) if and only if $T$ is also $\f$-positive.}

\pf   
Note that $T(\f)=\int\f(\oa T)d\|T\| \leq \int d \|T\| = M(T)$ since $\f(\oa T)\leq \|\oa T\|=1$.
Equality occurs if and only if $\f(\oa T)=1$  ($\|T\|$-a.e.).  This is Condition 3) in Lemma \DD.5.
The second assertion follows from the fact that $T_0(\f)=T(\f)$.  \qed

\medskip

 The reader may note that only Condition 3) of \DD.5 was used in this proof.
 However, it is Conditions 1) and 2) which
give a genuine understanding of $\f$-positive currents.

The closed currents which are $\f$-positive have a  {\sl monotonicity property}
which says that the function $\|T\|(B_x(r))/r^p$ is monotone increasing in $r$.
This implies that the density of $\|T\|$ is well-defined everywhere and upper
semi-continuous. This is disussed in detail in [HL$_3$].

Deep results in geometric measure theory
have important applications here.

\Theorem{\DD.8} {\sl Fix a compact set $K\subset X$ and a constant $c>0$. Then the set
$\cp(\f, K,c)$ of $\f$-positive currents $T$ with $M(T)\leq c$ and $\supp(T)\subseteq K$ is compact
in the weak topology.  }
\pf    Proposition  \DD.7 easily implies that a weak limit of $\f$-positive currents is $\f$-positive.
The result then follows from standard compactness theorems for measures.\qed

\Theorem{\DD.9} {\sl Fix a compact set $K\subset X$ and a constant $c>0$. Then the set
${\cal R}(\f, K,c)$ of  rectifiable $\f$-currents $T$ with rectifiable boundaries, such that $M(T) + M(\partial T) \leq c$ and $\supp(T)\subseteq K$, is compact
in the weak topology.}
\pf
This  follows from Proposition \DD.7 and the 
Federer-Fleming weak compactness theorem for rectifiable currents [FF], [F$_1$].
\qed

\Theorem{\DD.10} {\sl  Let $T$ be a  $\f$-cycle on $X$.  Then
there is a closed subset $\Sigma\subset \supp(T)$ of Hausdorff  dimension $p-2$ such that
$M \equiv \supp(T)-\Sigma$ is a proper $\f$-submanifold with finite volume in $X-\Sigma$ and 
$$
T\ =\ \sum_k n_k [M_k]
$$
where the $n_k$'s are positive integers and the $M_k$'s are the connected components of $M$.
}
\pf This is a direct consequence of Almgren's  big regularity theorem [A].\qed \medskip

We now present the dual characterization of $\f$-positive currents.

Suppose $\f\in \L^p V$ is a calibration on an inner product space $V$.
Let $\lpp\subset \Lambda^p V$ denote the polar cone of $\lp\subset \Lambda_p V$.
By definition this is the set of $\a\in \L^pV$ such that $\a(\x)\geq0$ for all $\x\in \L_+(\f)$, or 
equivalently,
$$
\lpp = \{\a \in \L^p V : \a(\x)\geq 0 {\ \ \rm for\ all\ } \x\in
G(\f)\}.
$$
A   $p$-form $\a\in \L^pV$  is said to be {\bf  $\lpp$-positive} 
if  $\a \in \lpp$, and {\bf  strictly $\lpp$-positive} if $\a(\x)>0$      
for all $\x\in G(\f)$ (or equivalently, $\a$ belongs to the interior of
 $\lpp$).

\Remark{\DD.11} Note that $\f$ itself is strictly $\lpp$-positive,
i.e.,  an interior point of the cone $\lpp\subset \L^p V$.
 If a closed convex cone has one interior point, then there
exists a basis for the vector space consisting of interior points.
Consequently, $\L^p V$ has a basis of strictly $\lpp$-positive $p$-forms.

\medskip

If $(X,\f)$ is a calibrated manifold, the considerations and definitions
above apply to the tangent space $V=T_x X$ at each point $x\in X$. 

\Def{\DD.12}   A smooth $p$-form $\a$ on $X$ is {\bf  $\lpp$-positive (strictly
$\lpp$-positive) } if $\a$ is  $\lpp$-positive (strictly $\lpp$-positive) at each point  $x\in X$.

\Def{\DD.13} A (twisted) current $T$ of dimension $p$ is said to be {\bf  $\lp$-positive}
if  $$T(\a)  \geq  0$$ for all $\lpp$-positive $p$-forms $\a$ with compact
support.

\Theorem{\DD.14} {\sl  A current $T$ is $\lp$-positive if and only if it is
 $\f$-positive.}

\medskip

This result is proven in [HL$_3$, Prop. A.2 and Remark on page 83]. 
However, for the sake of completeness we include a proof.

\pf First assume that $T$ is representable by integration.  Then 
from (2.1) we see that $T$ is $\L_+(\f)$-positive if and only if 
$$
T(g\a)\ =\ \int g \a(\oa T) d\|T\| \ \geq \ 0
$$
for all functions $g\geq 0$ and all compactly supported $\lpp$-positive $p$-forms $\a$.
Equivalently, each measure $\a(\oa T)\|T\|$ is $\geq 0$ for the same set of $p$-forms.  In turn, this is equivalent to
$$
\a(\oa T)\ \geq \ 0 \qquad \|T\| {\rm -a.e.}.
$$
for all compactly supported $\lpp$-positive $p$-forms.  Finally, by the Bipolar Theorem [S]
this last condition is equivalent to the Condition 1) in the Lemma \DD.5.

It remains to prove that if $T$ is $\L_+(\f)$-positive, then $T$ is representable by integration.
For this we may assume that $T$ has compact support
in a small neighborhood $U$ of $X$, and by    Remark \DD.11, we may choose a
frame $\a_1,...,\a_N$ of smooth $p$-forms which are strictly $\lpp$-positive
on $U$.  Let $\x_1,...,\x_N$ denote the dual frame of $p$-vector fields,
i.e., $(\a_i, \x_j)\equiv\d_{ij}$ on $U$.  Every such current $T$ has a
 unique representation as $T=\sum_{j=1}^N u_j \x_j$ with
$u_j \in \cd'(U)$ a distribution defined by $u_j(f) \equiv T(f\a_j)$ for all
test functions $f$. (Note that $\a=\sum_jf_j\a_j$ implies that $T(\a) =
\sum_jT(f_j\a_j) = \sum_j u_j(f_j) = (\sum_j u_j \x_j)(\sum_i f_i\a_i)
= (\sum_j u_j \x_j)(\a)$.)
Since $T$ is $\lp$-positive, each $u_j$ satisfies 
$$
u_j(f)\ \geq\ 0 \foral f\geq 0.
$$
By the Riesz Representation Theorem this proves that each $u_j$ is a
measure.  Therefore $T=\sum_j u_j \x_j$ is representable by integration.
\qed
\medskip

Now we give the proof of Lemma \DD.5.  As before $\f\in \L^pV$ is a calibration.
Let $K$ denote the unit mass ball in $\L_p V$, that is, the convex hull of 
the Grassmannian $G(p,V)\subset \L_p V$.

\Lemma {\DD.15}  $$ {\rm ch}\, G(\f) \ =\ \{\f=1\}\cap \partial K \ =\ \L_+(\f) \cap \partial K  $$

\pf  Note that:
\smallskip

a) \ \   $ {\rm ch}\, G(\f)   \subset   \{\f=1\}$\qquad since $G(\f)\subset \{\f=1\}$.

\smallskip

b) \ \   ${\rm ch}  \, G(\f)  \subset  K$\qquad\qquad \ \ since $G(\f)\subset G(p,V)$.

\smallskip

c)  \ \ $K\cap \{\f=1\}\ =\ \partial K\cap \{\f=1\}$ \ \ since $K\subset \{\f\leq 1\}$.

\smallskip\noindent
Hence, $ {\rm ch} \, G(\f)  \subset \{\f=1\}\cap\partial K$.

Conversely, suppose $\f(\x)=1$ and $\|\x\|=1$.  Since $\x\in K$,
$$
\x\ =\ \sum_j \lambda_j \x_j \qquad {\rm with\ each\ \ }\x_j \in G(p,V), \ {\rm each\ } \lambda_j >0, \ {\rm and\ }\sum_j \lambda_j=1.
$$
Hence, $1=\f(\x) = \sum \lambda_j\f(\x_j) \leq \sum\lambda_j=1$ forcing each $\f(\x_j)=1$ and therefore each $\x_j\in G(\f)$.  This proves the first equality.

We have shown ${\rm ch}  \, G(\f)   \subset \partial K$, and by definition, 
${\rm ch}  \, G(\f)  \subset \L_+(\f)$.  Hence, $ {\rm ch}\, G(\f)  \subset \lpp\cap \partial K$.
Finally, suppose $\x \in \partial K \cap \L_+(\f) $,  i.e., $\|\x\|=1$ and there exists   some $\lambda >0$ such that
$\lambda\x\in {\rm ch} \, G(\f)$.  We have already shown that 
 $ {\rm ch} \, G(\f)  \subset \partial K$, therefore $\|\lambda\x\|=1$, and hence 
 $\lambda =1$ proving that $\x\in{\rm ch} \, G(\f)$.
\qed

\Cor{\DD.16}  {\sl  Suppose $\x \in \L_pV$ has mass norm $\|\x\|=1$.  Then}
$$
\x\in \L_+(\f) \quad\iff\quad   \f(\x)=1 \quad\iff\quad  \x\in {\rm ch}  \, G(\f).
$$
This is the required restatement of Lemma \DD.5

\Remark{}  Note also  that the equation
$$
G(\f)\ =\ G(p,V)\cap \L_+(\f)
\eqno{(\DD.4)}
$$
follows easily from Lemma \DD.15.  This clarifies the notion of a rectifiable $\f$-current.
Namely, this proves that a rectifiable current is $\lp$-positive
if and only if it is a rectifiable $\f$-current, and eliminates a potential conflict in terminology.
\medskip

We finish this section with a lemma and corollary that are often useful.

Note that a  form $\a \in \Lambda^+(\phi)$ lies on the topological boundary $\partial \Lambda^+(\phi)$
 of $\lpp$ if and only if there
exists some  $\x\in G(\f)$ with $\a(\x)=0$.  

\Lemma{\DD.17} {\sl For any  $\psi\in \L^p V$}
$$
\f-\psi \in \partial  \lpp  \ \ \ \iff  \ \ \  \psi \leq 1 {\rm \ on\ } G(\phi) {\ \rm\ and \ }  \psi(\xi) =1 {\rm \ for\  some\   }  \xi \in G(\phi)
$$
\noindent
{\bf Proof.}  By definition $\phi-\psi \in \lpp$ if and only if $\phi(\xi)-\psi(\xi) = 1-\psi(\xi) \geq 0$ for all 
$\xi\in G(\phi)$. As remarked above $\phi-\psi$ lies in the boundary of $\lpp$ iff 
$\phi(\xi)-\psi(\xi) = 1-\psi(\xi) = 0$ at some point $\xi  \in G(\f)$.
\qed
\medskip

\Cor{\DD.18}  {\sl  For each unit vector $e\in V$, let $\f_e = e\hk(e\wedge\f) \cong \f\bigr|_W$ where
$W\equiv (\span e)^{\perp}$.  Then  $\phi \in \lpp$ and}
$$
\f_e \in \partial \lpp \qquad\iff \qquad  e\in \span \x \ \ {\rm for\ some\ }\ \x\in G(\f).
$$

\pf Note first that $\f_e = e\hk(e\wedge\f) = \f - e\wedge(e\hk \f)$.
Now write $e=a+b$ with $a\in \span\x$ and $b\perp \span\x$.
Then
$$
e\wedge(e\hk \f)(\x) = |a|^2.
\eqno{(\DD.5)}
$$
because
  $e\wedge(e\hk \f)(\x) = \f (e\wedge(e\hk \x)) = \f ((a+b)\wedge(a\hk \x)) = 
\f (|a|^2 \x) + \f( b\wedge(a \hk \x)) = |a|^2$ since 
$\f( b\wedge(a \hk \x)) =0$ by the First Cousin Principle [HL$_5$, Lemma 2.4] and $\f(\x)=1$.

Now by Lemma \DD.17,  $\f_e \in\partial \lpp $ if and only if  there exists $\x\in G(\f) $ with $|a|=1$, that is, with $e=a\in \span \x$.\qed

\Remark{}  Both $df\wedge(\nabla f\hk \f) = df\wedge d^\f f$  and  $\nabla f\hk (df \wedge \f) = \|\nabla f\|^2 \f - df \wedge d^\f f$ take values in $\lpp\subset \L^pT^*X$.  Furthermore, 
\smallskip

1) \ \ $df\wedge d^\f f \ \in \  \partial \lpp  \qquad \iff  \qquad \exists\ \x \in G(\f)$ tangential to the level set
 of $f$.
\smallskip

2) \ \  $\nabla f \hk (df\wedge \f)\ \in \  \partial \lpp 
  \qquad \iff  \qquad \exists\ \x \in G(\f)$ with $\nabla f \in \span \x$.
  \medskip
  \noindent
  Note that for some calibrations, condition 2) is true for all $f$, i.e., given a vector $n\in V$, there always exists a $p$-vector $\x \in G(\f)$ with $n\in \span \x$.

\vskip .3in

\centerline{\bf \EE.   Duality with Plurisubharmonic Functions.}\medskip

As noted in the introduction,  the $\f$-plurisubharmonic functions represent,  in a general sense,  the polar dual of the $\f$-positive currents.  More specifically, there are many    situations where
the polar dual of an interesting  family of $\f$-positive currents is some explicit and useful family of 
$\f$-plurisubharmonic  functions. In this section we examine an example of this phenomenon.
We extend the fundamental duality results of Duval-Sibony  [DS]  in complex geometry, to a general calibrated manifold  $(X,\f)$. The  Duval-Sibony Duality Theorem involves plurisubharmonic functions,
pseudoconvex hulls, positive currents and Poisson-Jensen formulas.

\Def{\EE.1}  Suppose $R$ is a $(p-1)$-dimensional current on $X$.  The operator $\bdf$ is defined by
$$
(\bdf R)(f)\ =\ R(d^\f f)
$$
for all $f\in C^{\infty}_{\rm cpt}(X)$.

In other words, $\bdf : \ce_{p-1}'(X) \arr \ce_0'(X)$ is the adjoint of $d^\f:\ce^0(X) \arr \ce^{p-1}(X)$.
Let $\bd :\ce_{p}'(X)\arr \ce_{p-1}'(X)$ denote the {\bf boundary operator} on currents.  This is the  
 adjoint of $d:\ce^{p-1}(X) \arr \ce^p(X)$ and is related to the deRham differential on currents by
$\bd =(-1)^{n-p}d$.  The adjoint of $dd^\f:\ce^0(X)\arr \ce^p(X)$ is the operator
$$
\bdf \bd : \ce_p'(X) \ \arr\ \ce_0'(X)
\eqno{(\EE.1)}
$$

\Remark{}   Throughout the remainder of this section we assume that $(X,\f)$ is a non-compact connected calibrated manifold.  We also assume that $\f$ is parallel.  This assumption enables
us to use the operator $\partial_\f \partial$, but it is not necessary.  We leave it to the reader
to  verify that all of the results of this section extend to the case where $\f$ is not parallel  by replacing
 the operator  $\partial_\f \partial$ with $\ch_\f : \ce_p'(X)\to  \ce_0'(X)$, the  adjoint of 
 $\ch^\f:\ce^0(X)\to \ce^p(X)$.  Of course, $\ch_\f$ is defined by 
 $(\ch_\f(T))(f) = T(\ch^\f(f))$ for all $f\in C^{\infty}_{\rm cpt}(X)$.

\Lemma{\EE.2. (The Support Lemma)}  {\sl  Suppose $K$ is a compact subset of $X$.  Suppose $T$ is a $\lp$-positive current with compact support in $X$.  If $\bdf \bd T$ is $\leq 0$ (a non-positive measure) on $X-\wh K$, then 
$\supp \,T \subset \wh K\cup \Core(X)$.}

\pf  The Note following Lemma 4.2 in [HL$_5$] states that for each $x\notin \wh K\cup \Core(X)$ there exists a non-negative \fp function $f$ on $X$ which is identically zero on a neighborhood of $\wh K$  and strict at $x$. 
 Since $f$ is strict  at $x$, there exists a small ball $B$ about $x$ and $\e>0$
  so that $dd^\f f-\e \f$ is $\lpp$-positive at each point in $B$. 
   By equation (\DD.2), $M(\chi_B T)= (\chi_B T, \f)$.  Therefore, 
$\e M(\chi_B T) = (\chi_B T, \e \f) \leq (\chi_BT , dd^\f f)\leq (T,dd^\f f) = (\bdf \bd T, f)\leq 0$.
This proves that $M(\chi_B T)=0$  and hence $\supp\, T\subset \wh K\cup \Core(X)$.\qed
\medskip

The case where $K=\emptyset$ is a generalization of  Proposition 4.13 in [HL$_5$] .
\Cor{\EE.3} {\sl  If $T$ is a $\lp$-positive current with compact support and $\bdf\bd T\leq 0$, then}
$$
\supp\, T\ \subset \  \Core(X).
$$

When $\Core(X)=\emptyset$ we have

\Cor{\EE.4} {\sl  Suppose $(X,\f)$ is strictly convex and $K$ is $\f$-convex.  
Suppose $T$ is $\lp$-positive with compact support.  If $\supp \{\bdf\bd T\}\subset K$, then $\supp\, T\subset K$.
In particular, there are no $\lp$-positive currents which are compactly supported without boundary on  $X$.}
\medskip
.

We now introduce the notion of a Green's current on $(X,\f)$.
We shall begin with a description of the ``classical'' case.
Let $M\subset X$ be an $p$-dimensional oriented submanifold  having no compact components, 
and fix a compact domain $D\subset\subset M$ with smooth boundary $\partial D$.
(One can assume that M is just $D$ with an external collar added to it.)
 Let $G_x$ denote the Green's function for $(D, \partial D)$ with singularity at $x\in {\rm Int} D$. 
 Let  $\mu_x$ denote harmonic measure (i.e., the Poisson kernel)  on $\partial D$ and let $[x]$ denote the  point-mass measure at $x$.  Extend $G_x$ to a continuous function on $M$ (also denoted
 by $G_x$) by defining it to be zero on $M-D$.   Then
$$
*_M \Delta_M G_x\ =\ \mu_x - [x]\qquad {\rm on\ \ }  M.
\eqno(\EE.2)
$$

If $M$ is a $\f$-submanifold of the  calibrated manifold $(X,\f)$, this  equation can be reformulated as a
current $\bdf \bd$-equation on $X$.

\Lemma{\EE.5}  {\sl  Suppose $M$ is a $\f$-submanifold of $X$, and that $u\in
\cd'^0_{\rm cpt}(M)$ is a generalized function with compact support on $M$.  Then  }
$$
\bdf \bd(u[M])\ =\ (*_M\Delta_M u)[M].
$$

\pf
Consider the inclusion map $i:M\hookrightarrow X$.  Then, by definition,
$u[M] = i_* u$ and $(*_M\Delta_M u)[M] = i_*(*_M\Delta_M u)$.  For any test
function $f$ on $X$ we have 
$  
((\bdf \bd)(i_* u), f) \ =\ (i_*u, dd^\f f)\ =\ (u, i^*(dd^\f f))_M
$
where $(\cdot , \cdot)_M$ denotes the pairing of functions with currents on
$M$.   Equation (1.2) of [HL$_5$]  states that $i^*(dd^\f f) \ =\ *_M \Delta_M(i^* f)$
if $M$ is a $\f$-submanifold.   Finally,
$
(u,  *_M \Delta_M(i^* f))_M\ =\ ( *_M \Delta_M u, i^* f)_M
\ =\ (i_*( *_M \Delta_M u), f)_X.
$
\qed

\Cor{\EE.6} {\sl  Suppose $M$ is a
$\f$-submanifold and $G_x$ is defined as above.  Then
$$
\bdf \bd (G_x[M])\ =\ \mu_x - [x]
$$
as a current equation on $X$.}
\medskip

We now generalize these currents $G_x[M]$.
Assume $K$ is a compact subset of $X$, and let $\cm_K$ denote the set of
probability measures with support in $K$.  

\Def{\EE.7}  If $T_x$ is a $\lp$-positive current with compact support and
$T_x$ satisfies:
$$
\bdf \bd \,  T_x\ =\ \mu_x -[x]
\eqno(\EE.3)
$$
with $\mu_x\in \cm_K$, then: $T_x$ is a {\bf Green's current  for }
$(K,x)$, $\mu_x$ is a {\bf Poisson-Jensen measure} for $(K,x)$, and the
equation (\EE.3) is the {\bf Poisson-Jensen equation}.

\medskip

\Theorem{\EE.8} {\sl  Suppose $X$ is strictly $\f$-convex, $K$ is a
compact subset of $X$, and $x\in X-K$. Then there exists a Green's current
$T_x$ for $(K,x)$ if and only if $x\in \PH K$.
}

\medskip
To prove this we begin with the following.

\Prop{\EE.9}  {\sl  Suppose $(X,\f)$ is noncompact   calibrated manifold. 
If there exists a Green's current for $(K,x)$, then $x\in \PH K$.
}

\pf 
Since $\bdf \bd \, T_x \ =\ \mu_x-[x]$, we have $\int f \mu_x - f(x) \ =\
(T_x, d d^\f f)$ for all $f\in C^\infty(X)$.
If $f$ is \fp, this implies that $f(x)\ \leq\  \int f \mu_x\ \leq \ \sup_K
f$, since $\mu_x\in \cm_K$. Thus $x\in \PH K$.\qed
\medskip

The set $\cp_X \equiv \fpsh\subset C^\infty(X)$  of all
\fp  functions on $X$ is clearly a closed convex cone in
$C^\infty(X)$.  Let 
$$
 C_X\ \equiv\ \{u\in \ce'_0(X) : u=\bdf \bd \, T \ {\rm for\ some\ } \lp{\rm-positive\ }T\in
\ce'_p(X)  \}. 
\eqno{(\EE.4)}
$$
This is a convex cone in $ \cd'_{0, {\rm cpt}}(X)$.
Let's abbreviate $\cp = \cp_X$ and $C=C_X$.

\Lemma{\EE.10} {\sl  Suppose $X$ is non-compact with calibration $\f$.
Then $\cp$ is the polar of $C$,
that is,
$$
\cp \ =\ C^0\  \equiv \ \{f\in C^\infty(X) \ : \  (u,f)\geq 0 \ \  \forall u
\in C\}. $$
}
\pf
Consider $u=\bdf \bd(\delta_x \x)$, with $\x\in G_x(\f)$.  Clearly $u\in C$. 
If $f\in C^\infty(X)$ belongs to $C^0$, then 
$0\leq (u,f) = (\bdf \bd (\delta_x \x), f) \ =\ (\delta_x \x, d d^\f f)\ =\
(dd^\f f)_x(\x).$  Hence $C^0\subseteq \cp$.

Conversely, if $f\in\cp$, then for all $u\in C$, $(u,f) \ =\ (\bdf \bd   T, f)\ =\ (T, dd^\f
f)\geq 0$, since $T$ is $\lp$-positive. This proves that $\cp\subseteq C^0$.
\qed

\Lemma{\EE.11}  {\sl  If $X$ is strictly $\f$-convex, then the convex cone
$C\subset \cd'_{0, {\rm cpt}}(X)$ is closed.} 

\pf
Let $\ce_{0,K}'(X) = \{T\in \ce_0' (X):  {\rm supp}\, T\subset K\}$.
It suffices to show that $C\cap \ce'_{0, K}(X)$ is closed for an exhaustive
family of compact subsets $K\subset X$.  We may assume $K$ is $\f$-convex.
Suppose $u_j$ converges to $u$ in $C\cap \ce'_{0, K}(X)$ with each $u_j\in C$,
i.e., $u_j = \bdf \bd T_j$ where $T_j$ is a $\lp$-positive current with
compact support.  By Corollary \EE.4  the support of each $T_j$ is
contained in $K$.  Consider a strictly \fp function $f$ on $X$.  Pick
$\epsilon>0$ so that $dd^\f f -\epsilon \f$ is $\lpp$-positive at each
point of $K$.  Then 
$
\epsilon M(T_j) \ =\ (T_j, \epsilon \f) \ \leq\ (T_j, dd^\f f)\ =\ (\bdf \bd 
T_j, f)\ =\ (u_j,f) 
$
which converges to $(u,f)$.  Therefore the masses $M(T_j)$ are bounded.
By compactness there exists a weakly convergent subsequence $T_j\to T$.
Now $\supp\, T\subset K$ and $T$ must be $\lp$-positive. 
Hence $u=\bdf \bd T\in C\cap \ce'_{0, K}(X)$.  
This proves that $C\cap \ce'_{0, K}(X)$ is closed for each compact set $K$
which is $\f$-convex. \qed

\Cor{\EE.12}  {\sl  Suppose $X$ is strictly $\f$-convex.  Then
$$
C \ =\ \cp^0.
$$
Equivalently, the equation 
$$
\bdf \bd  T\ =\ u
$$
 has a solution $T$ which is a
$\lp$-positive current with compact support if and only if}
$$
0\ \leq \  u(f)\qquad {\sl for\ all\ \ \ } f\in\fpsh
$$
\pf 
Apply the Bipolar Theorem (cf. [S]).\qed

\medskip
\noindent
{\bf Proof of Theorem \EE.8.}
Suppose there does not exist a Green's current for $(K,x)$, that is,
suppose $\cm_K-[x]$ is disjoint from the cone $C$.  By the Hahn-Banach
Theorem (note that $\cm_K-[x]$ is a compact convex set) there exists  $f\in
C^0=\cp$ with $f$, considered as a linear functional on $\cd'_{0, {\rm cpt}}(X)$, 
satisfying $(u,f) \leq -\e <0$ for all $u\in (\cm_K-[x])$.   
That is, $\int f d\,\mu - f(x) \leq -\epsilon <0$ for all $\mu\in
\cm_K$. Consequently, 
$$
\sup_K f \ =\ \sup_{\mu\in \cm_K} \int f  d\,\mu\ \leq\ f(x) -\epsilon
$$
or $\sup_K f +\epsilon \leq f(x)$ so that $x\notin \PH K$.
In light of Proposition \EE.9 we are done. \qed
\medskip

One could define the ``Poisson-Jensen hull'' of a compact set $K$ to be the set of points $x$
for which there exists a Poisson-Jensen measure $\mu_x$ and a Green's current $T_x$ satisfying
(\EE.3).
Then Proposition \EE.9 states that on any (non-compact) calibrated manifold $(X,\f)$, the 
Poisson-Jensen hull of a compact set is contained in the \fp hull, while Theorem \EE.8
states that the two hulls are equal if $(X,\f)$ is strictly convex.

The next ``hull"  obviously contains the Poisson-Jensen hull.

\Def{\EE.13}  The {\bf current hull} of a compact subset $K\subset X$ is the union
$$
\wt K\ \equiv\ \bigcup_{T\in \cp(K)}
\supp\, T
$$
where $\cp(K)$ consists of all $\lp$-positive currents with compact support on $X$ 
satisfying $\bdf \bd  T\leq 0$ on $X-K$.

\Lemma{\EE.14}  {\sl  If $(X,\f)$ is strictly $\f$-convex, then $\wt K = \wh K$.}

\pf
The Support Lemma \EE.2 states that 
$$
\wt K\ \subset \ \wh K \cup \Core(X).
$$
for any calibrated manifold.  Now $\wt K$ contains the Poisson-Jensen hull
which equals $\wh K$ if $(X,\f)$ is strictly \fc  by Theorem \EE.8.
Strict $\f$-convexity also implies $\Core(X)=\emptyset$ 
by Theorem 1.3.  \qed

Suppose $(X,\f)$ is non-compact and connected.

\Def{\EE.15}
  An open subset $\Omega \subset X$ is   {\bf \fc  relative to $X$} if $K\subset\subset \Omega$
  implies ${\wh K}_X\subset\subset \Omega$.\medskip
  
  Note that if $X$  is \fc this condition implies that $\Omega$ is \fc since
 ${\wh K}_\O\subseteq{\wh K}_X$.
 Moreover, if $X$ is strictly \fc, then $\Core(\O)\subseteq\Core(X)$ is empty so that $\O$ is
 strictly \fc (by Theorem 1.3).

\Prop{\EE.16}  {\sl  Suppose $(X,\f)$ is strictly \fc.  An open subset  $\O \subset X$ is \fc relative to $X$ if and only if  ${\rm PSH}(X,\f)$ is dense in ${\rm PSH}(\O,\f)$.}

\pf Let $L:\ce^0(X)\to \ce^0(\O)$ denote restriction. The adjoint $L^*:\ce_0'(\O) \to 
\ce_0'(X)$ is inclusion.  Suppose $v\in (L^*)^{-1}(C_X)$,  i.e., $v\in \ce_0'(\O)$ with $v=\partial_\f\partial T$ for some $\lp$-positive current $T$ compactly supported in $X$.  
Then $K\equiv \supp\, v \subset \Omega$ satisfies $\wt K_X =\wh K_X$ by Lemma \EE.14. 
 Hence, $\wh K_X\subset \O$ implies    $\supp\, T \subset \O$, i.e.,  that $v\in C_\O$.  
 This proves that $\O$ is \fc relative to $X$ if and only if 
$$
 (L^*)^{-1}(C_X) \ =\  C_\O.
\eqno{(\EE.5)}
$$
By Corollary \EE.12 we may replace $C_X$ by $\cp_X^0 = C_X$.  In general, 
$[L(\cp_X)]^0 = (L^*)^{-1}(\cp_X^0)$, so that (\EE.5) is equivalent to
$$
\left[ L(\cp_X) \right]^0 \ =\ C_\O.
\eqno{(\EE.6)}
$$
By Lemma \EE.10, $\cp_\O = C_\O^0$.  Hence (\EE.6) is equivalent to 
$$
\overline{L(\cp_X)} \ =\ \cp_\O
\eqno{(\EE.7)}
$$
\qed



\vskip .3in

\centerline{\bf \HH.  \  Boundary Duality.  }

\medskip

In this section we take up the following general question.  Suppose $(X,\f)$ is a non-compact strictly \fc manifold.  Given a compact oriented submanifold $\G\subset X$ of dimension $p-1$, when does there exist a $\f$-submanifold $M$ with boundary $\G$?  More generally, when does there exist
a  $\lp$-positive current $T$ with $\partial T=\G$?

\Theorem{\HH.1} {\sl Suppose $\f$ is exact.  Given $S\in \ce'_{p-1}(X)$, there  exists a $\lp$-positive current $T\in \ce'_p(X)$ with $S=\partial T$ if and only if }
$$
\int_S\a \ \geq\ 0\qquad {\sl for\ all\ }  \a \in \ce^{p-1}(X) \ {\sl such \  that \ } d\a \ {\rm is\ } \ \lpp{\sl -positive}
$$

\pf
Consider the following convex cones.
$$
A\ =\ \{\a \in \ce^{p-1}(X) : \  d\a \ {\rm is\ } \ \lpp{\rm -positive}\}
$$
$$
B\ =\ \{S\in \ce'_{p-1}(X) : \ S=\partial T \ {\rm for\ some\ } \lp{\rm -positive} \ T\in \ce'_p(X)\}
$$

If $\a\in A$ and $S\in B$, then
$$
S(\a)=\partial T(\a) = T(d\a)\geq 0,
$$
that is,
$$
A\ \subseteq \ B^0 \and B\ \subseteq \ A^0
$$
where $B^0$ denotes the polar of $B$.  If $\x\in G_x(\f)$, then $T=\delta_x \x$ is $\lp$-positive, so that $\partial(\delta_x \x)\in B$. Therefore, if $\a\in B^0$, then $0\leq \partial(\delta_x \x)(\a)
=(\delta_x \x)(d\a) = (d\a)_x(\x)$.  This proves that $B^0\subseteq A$, and hence $A=B^0$.
(In particular, note that $A$ is closed.)  Theorem \HH.1 is just the statement that
$B=A^0$.  Now since $A=B^0$, the Bipolar Theorem states that $\overline B = A^0$.
Thus it remains to show that $B$ is closed.

Suppose $S_j\in B$ and $S_j\to S$ in $\ce'_{p-1}(X)$.  Then $S_j=\partial T_j$ for some 
$T_j$ which is $\lp$-positive.  The calibration $\f$ is exact, i.e., $\f=d\eta$ for some
 $\eta\in\ce^{p-1}(X)$.  Therefore
 $$
 M(T_j)\ =\ T_j(\f)\ =\ T_j(d\eta)\ =\ (\partial T_j)(\eta)\ = \ S_j(\eta)\ \arr \ S(\eta).
 $$
 In particular, there exists a constant $C$ such that $M(T_j)\leq C$ for all $j$.
 By the Support Lemma \EE.2 and Theorem 1.3 we have
 $$
 \supp T_j\ \subseteq\ \wh{\supp S_j}
 $$
for each $j$.  Pick a compact subset $K$ with $\supp S_j \subseteq K$ for all $j$.
Then
$$
 \supp T_j\ \subseteq\ \wh{K}\qquad {\rm for\ all\ } j.
$$

This proves that $\{T_j\}$ is a precompact set in $\ce'_p(X)$.  Therefore, there exists a convergent
subsequence $T_j\to T$ in $\ce'_p(X)$. Obviously, $\partial T=S$ and $T$ is $\lp$-positive.
Hence, $S\in B$.
\qed

\Remark{\HH.2}  The same proof combined with the compactness Theorem \DD.9 
 proves the following.  Let $\ccr_p(X)$ denote the compactly supported rectifiable 
 currents of dimension $p$ on $X$. Then, if $\f$ is exact, the set
$$
B_{\rm rect} \
 \equiv\ \{ \G\in \ccr_{p-1}(X) : \ S=\partial T \ {\rm for\ some\ } \lp{\rm -positive} \ T\in \ccr_p(X)\}\}
$$
is weakly closed in $\ccr_{p-1}(X)$.

\medskip 

We now turn attention to the case where $X$ is compact.
\def\C{\lambda}

\Theorem{\HH.3} {\sl Suppose $(X, \f)$ is a compact calibrated manifold.  Fix  $S\in \ce'_{p-1}(X)$
and $\C>0$.    Then  the following are equivalent.  
\smallskip
(i) \ There  exists a $\lp$-positive current $T\in \ce'_p(X)$ with $S=\partial T$  and    
$M(T)\leq \C$.  \smallskip

(ii)\  $\int_S\a \ \geq\ -\C$ for all $\a \in \ce^{p-1}(X)$ such that $d\a+\f$ is  $\lpp$-positive.}
 
\pf  Consider the convex sets:
\smallskip
\centerline{$A_{\C}\ =\ \{\a\in  \ce^{p-1}(X) : \  d\a+\smfrac 1 \C \f$  \ \  is  $\lpp$-positive$\}$,}
\medskip

\centerline{
$B_{\C}\ =\  \{S\in \ce'_{p-1}(X) : \ S=\partial T$
 \ \  for some $\lp$-positive  $T\in \ce'_p(X)$
with  $M(T)\leq \C\}$.
}
\medskip\noindent
If $\a\in A_{\C}$ and $S\in B_{\C}$, then
$$
S(\a)=\partial T(\a) = T(d\a) = T(d\a +\smfrac 1 \C \f -\smfrac 1 \C \f )\geq 
-\smfrac 1 \C T( \f ) = -\smfrac 1 \C M(T) \geq -1,
$$
that is,
$$
A_{\C}\ \subseteq \ B_{\C}^0 \and B_{\C}\ \subseteq \ A_{\C}^0
$$
where $B_{\C}^0=\{\a: (S,\a)\geq -1$ for all $S\in B_{\C}\}$ 
 denotes the polar of the convex set  $B_{\C}$.  
 If $\x\in G_x(\f)$, then $T={\C}\delta_x \x$ is $\lp$-positive, so that 
 $\partial({ \C}\delta_x \x)\in B_{\C}$. Therefore, if $\a\in B_{\C}^0$, 
 then $-1  \leq \partial(\lambda\delta_x \x)(\a)
= \lambda (\delta_x \x)(d\a) = \lambda(d\a)_x(\x)$.   Since $\f(\x)=1$ we conclude that
$(d\a+{1\over \C}\f)(\x)\geq0$.
This proves that $B_{\C}^0\subseteq A_{\C}$, and hence $A_{\C}=B_{\C}^0$.
(Note, in particular,  that $A_{\C}$ is closed.)  Theorem \HH.3 is just the statement that
$B_{\C}=A_{\C}^0$.  Now  the Bipolar Theorem 
states that $\overline B_{\C} = A_{\C}^0$.
However it follows from the Compactness Theorem \DD.8 that 
$\overline B_{\C} = B_{\C}$.\qed

\medskip

Similar arguments can be applied to prove other versions of boundary duality.
For example we have the following result concerning boundaries of 
$\f$-positive currents whose support is not necessarily compact.

\Theorem{\HH.4} {\sl Suppose $(X, \f)$ is  an arbitrary calibrated manifold.  Fix  $S\in \cd'_{p-1}(X)$
and $\C>0$.    Then  the following are equivalent.  
\smallskip
(i) \ There  exists a $\lp$-positive current $T\in \cd'_p(X)$ with $S=\partial T$  and    
$M(T)\leq \C$.  \smallskip

(ii)\  $\int_S\a \ \geq\ -\C$ \ \  for all  compactly supported forms $\a \in \cd^{p-1}(X)$ such that $d\a+\f$    

\qquad\qquad\qquad\qquad is $\lpp$-positive.}

\medskip
Note that here the currents $S$ and $T$  do not necessarily have compact support.
The proof of Theorem \HH.4 is left to the reader.  Arguments for the special case of 
complex geometry appear in [HL$_4$].


\vfill \eject

\centerline{\bf \II.  \  $\f$-Flat Hypersurfaces and Functions which are $\f$-Pluriharmonic mod $d$. }

\medskip

The $\f$-pluriharmonic functions are the closest thing to holomorphic functions on a calibrated manifold $(X,\f)$.  Usually there are very few $\f$-pluriharmonic functions. An attempt has been made in this paper to remedy this situation by emphasizing the \fp functions. By comparison these functions exist in abundance.  For some purposes another extension of the concept of 
$\f$-pluriharmonic functions is more useful --- namely the {\sl $\f$-pluriharmonic functions mod $d$}.

This section is, for the most part, a straightforward extension of the results of Lei Fu [Fu]
from the special Lagrangian case to the general calibrated manifold $(X,\f)$.

\Def{\II.1}  A function $f\in C^\infty(X)$ is {\bf \fh mod $d$} if 
$$
dd^\f f \ =\ df\wedge\a_f + \s_f
\eqno{(\II.1)}
$$
for some $(p-1)$-form $\a_f$ and some $p$-form $\s_f$ of type $\L(\f)^\perp$, i.e.,
$\s_f(\x)=0$ for all $\x\in G(\f)$.
\medskip

If $f$ is \fh mod $d$, then $\lambda f$, $\lambda\in\bbr$, is also \fh mod $d$.  However, the sum of two such functions need not be \fh mod $d$.

\Prop{\II.2}  {\sl   Suppose that $df$ never vanishes so that $\ch \equiv \ker df$ defines a hypersurface foliation.  The condition that $f$ be \fh mod $d$ is independent of the function defining 
the foliation $\ch$.}

\pf
Recall that locally $f$ and $g$ define the same foliation $\ch$ if and only if $g = \chi(f)$ for some
function $\chi:\bbr\to\bbr$ for which   $\chi'$  is never zero.  To prove this fact
assume that $f=x_1$ is a local coordinate.  Since $g$ is constant on the leaves 
$\{x_1={\rm constant}\}$,  $g$ must be independent of $x_2,...,x_n$, i.e., $g=\chi(x_1)$.
Since $dg$ is never zero, $\chi'$ is never zero.  Finally, 
$$
dd^\f g\ =\ \chi'(f)dd^\f f +\chi''(f) df\wedge d^\f f
\ =\ dg\wedge\left(\a_f+{{\chi''(f) }\over{\chi'(f) }}d^\f f\right) +\chi'(f) \sigma_f
$$
which proves that if $f$ is \fh mod $d$, then $g=\chi(f)$ is also. \qed
\medskip

Recall that a real hypersurface $Y\subset X$ is {\sl $\f$-flat} if its second fundamental form
$II_Y$ has that property that $\tr_\x II_Y =0$ for all tangential $\f$-planes $\x\subset TY$.
This is equivalent to the condition that $d d^\f f(\x)=0$ for all such $\x$, where $f$ is any defining
function for $Y$.
(See Definition 5.1 and Lemmas 5.2 and 5.11 in [HL$_5$].)

\medskip\Prop{\II.3}  {\sl    If $f$ is \fh mod $d$, then each (non-critical) hypersurface $\{f=C\}$ is \ffl.}

\pf Suppose $\x\in G(\f)$ is tangent to $\{f=C\}$, i.e., $\nabla f \hk \x=0$.  Then
$
(dd^\f f)(\x) \ =\ (df\wedge \a_f)(\x) + \s_f(\x).
$
Since $\s_f$ is of type $\L(\f)^\perp$, we have $\s_f(\x)=0$, and $(df\wedge\a_f)(\x)= \a_f(\nabla f\hk \x)=0$. \qed\medskip

By  Corollary 2.11   in [HL$_5$] we have that for any $f\in C^\infty$ and any $\f$-submanifold $M$,
$$
(dd^\f f-df\wedge \a_f)\bigr|_M\ =\ *_M(\D_M f) - d(f\bigr|_M)\wedge\a_f \bigr|_M.
\eqno{(\II.2)}
$$
This proves the following.

\Prop{\II.4}  {\sl  If $f$ is \fh mod $d$ and $M$ is a $\f$-submanifold, then
$u\equiv f\bigr|_M$ satisfies the partial differential equation
$$
\D_M u \ =\ *(du\wedge\beta)   \qquad {\rm on\ \ } M
\eqno{(\II.3)}
$$
where $\b=\a_f\bigr|_M$.}\medskip

The maximum principle is applicable to solutions to (\II.3).  See  [BJS].

\Cor{\II.5}  {\sl  Suppose $(M,\G)$ is a compact $\f$-submanifold with boundary.
Then for each function $f$ which is \fh mod $d$ and each point $x\in M$, one has}
$$
\inf_\G f \ \leq \ f(x) \ \leq \ \sup_\G f
\eqno{(\II.4)}
$$

\Cor{\II.6}  {\sl   Suppose $(M,\G)$ is as above.  
If $\G\subset \{f=C\}$, then $M\subset \{f=C\}$.}

\Prop{\II.7}  {\sl  Suppose $(M,\G)$ is a compact $\f$-submanifold with boundary, and suppose
$f$ is a function on $X$ which is \fh mod $d$. If $f$ is constant on $\G$, then}
$$
d^\f f\bigr|_\G \ \equiv\ 0 \qquad{\rm (point-wise).}
\eqno{(\II.5)}
$$

\pf
By Corollary \II.6, $f$ is constant on $M$. We then apply the following.

\Lemma{\II.8}  {\sl For any function $f$ constant on $M$, $d^\f f\bigr|_\G \ \equiv\ 0$.}

\pf
At $x\in\G$, we have $\overrightarrow M = e\wedge \overrightarrow \G$ for some
$e$ tangent to $M$.  Since $f$ is constant on $M$, $\nabla f\perp \span \overrightarrow M$.
Now $(d^\f f)(\overrightarrow\G)= (\nabla f \hk \f)(e\hk \overrightarrow M) = \f((\nabla f)\wedge (e\hk\overrightarrow M))=0$ since $\nabla f \wedge (e\hk  \overrightarrow M)$ is a first cousin of 
$\overrightarrow M\in G(\f)$ (cf.  [HL$_5$, 2.4]). \qedqed

\medskip

Our next objective is to show that, for the large class of  {\sl normal} calibrations, a function
$f$ is \fh mod $d$ if and only if its level sets are \ffl.

Suppose $\f\in \L^pV$ is a calibration on a euclidean vector space $V$.  For each hyperplane
$W\subset V$, $\f\bigr|_W\in\L^p W$ has comass $\leq 1$ and, in fact, $<1$ if and only if
$G(\f\bigr|_W)$ is empty.

\Def{\II.9}  The calibration $\f\in \L^pV$ is {\bf normal} if, for every hyperplane $W\subset V$
$$
\L(\f\bigr|_W)^\perp\ = \L(\f)^\perp\bigr|_W
$$
as subspaces of  $\L^pW$.  A calibration $\f$ on a manifold $X$ is {\bf normal}
if $\f_x\in\L^p T_xX$ is normal for each $x\in X$.

\Prop{\II.10}  {\sl  Suppose $\f$ is a normal calibration on $X$, and $f\in C^\infty(X)$ has a         never-vanishing gradient.  Then
\smallskip

\centerline{\sl  $f$ is \fh mod $d$ }
\smallskip
\noindent
if and only if }

\smallskip

\centerline{\sl  each level set $\{f=C\}$  is \ffl.}
\smallskip

\pf Suppose each  level set $\{f=C\}$  is \ffl.  That is
$$
(dd^\f f)(\x) \ =\ 0 \qquad {\rm for\ all\ }\x\in G(\f) \ \ {\rm which\  are\ tangential\ to\ }
\{f=C\}.
\eqno{(\II.6)}
$$
Let $W=\ker df_x$ at some $x\in X$.
 Note that $G(\f\bigr|_W) = \{\x\in G(\f) : \x $ is tangential to $W\}$.   
Then (\II.6) is equivalent to 
$$
dd^\f f\bigr|_W \in \L\left(\f\bigr|_W\right)^\perp
\eqno{(\II.7)}
$$
Now $f$ is \fh mod $d$ if 
$$
dd^\f f\ =\ df\wedge \a_f +\s_f  \qquad \s_f \in \L(\f)^\perp
\eqno{(\II.8)}
$$
or equivalently
$$
dd^\f f\bigr|_W \in \L(\f)^\perp \bigr |_W
\eqno{(\II.9)}
$$
If $\f$ is normal, then $$ \L\left(\f\bigr|_W\right)^\perp  \subset  \L(\f)^\perp \bigr |_W$$
and (\II.7) implies (\II.9).
\qed

\Prop{\II.11}  {\sl The following calibrations are normal.

\medskip
 
\qquad 1.  A K\"ahler or  $p$th power K\"ahler calibration.

\medskip
 
\qquad 2.  A Special Lagrangian calibration.

\medskip
 
\qquad 3.  An associative, coassociative or Cayley calibration.

\medskip
 
\qquad 4.  A quaternionic calibration.

}
\medskip
\noindent
The proof is left to the reader.


\vfill\eject


\centerline{\bf 6.  \   Hodge Manifolds}

\medskip

\def\IH{\widetilde{H}}

\def\wt{\widetilde}
\def\smfrac#1#2{\hbox{${#1\over #2}$}}

\def\bbr{{\bf R}}
\def\bbh{{\bf H}}

\def\bbz{{\bf Z}}
\def\L{\Lambda}
\def\ce{{\cal E}} 
\def\cd{{\cal D}}

\font\AAA=cmr14 at 12pt
\font\BBB=cmr14 at 8pt

\def\O{\Omega}

\def\Theorem#1{\medskip\noindent {\AAA T\BBB HEOREM \rm #1.}}
\def\Prop#1{\medskip\noindent {\AAA P\BBB ROPOSITION \rm  #1.}}
\def\Cor#1{\medskip\noindent {\AAA C\BBB OROLLARY \rm #1.}}
\def\Lemma#1{\medskip\noindent {\AAA L\BBB EMMA \rm  #1.}}
\def\Remark#1{\medskip\noindent {\AAA R\BBB EMARK \rm  #1.}}
\def\Note#1{\medskip\noindent {\AAA N\BBB OTE \rm  #1.}}
\def\Def#1{\medskip\noindent {\AAA D\BBB EFINITION \rm  #1.}}

\def\Ex#1{\medskip\noindent {\AAA E\BBB XAMPLE \rm    #1.}}

\def\pf{\medskip\noindent {\bf Proof.}\ }

\def\qed{\hfill  $\vrule width5pt height5pt depth0pt$}

\def\lp{\Lambda_+(\phi)}
\def\lpp{\Lambda^+(\phi)}

In this section we pose some highly speculative questions for calibrated manifolds
in the spirit of those posed  in the complex case (cf. [HK, p.58], [L$_{4,5}$]).
Assume that $(X,\phi)$ is a compact calibrated $n$-manifold with a 
parallel calibration $\phi$ of degree $p$.  Let $*\phi$
denote the dual calibration.  Note that a $\phi $-submanifold or, more generally, any $\phi$-cycle
 (see Definition 2.1)
is a current of dimension $p$ and degree $n-p$.  
By contrast a $*\phi $-submanifold or $*\phi $-cycle is a current of 
dimension $n-p$ and degree $p$. 

 Denote by ${\IH}^p(X,\bbz)$  the image of  the map
${H}^p(X,\bbz) \to {H}^p(X,\bbr)$, with analogous  notation for homology.

\Def{6.1} If the de Rham class of the calibration $\phi$ lies in ${\IH}^p(X,\bbz)$, 
i.e., if $\phi$ has integral periods, then $(X,\phi)$ will be referred to as a {\bf Hodge manifold}.

\Remark{6.2}  If $(X,\omega)$ is a K\"ahler manifold, then this coincides with standard
terminology.  The Kodaira Embedding Theorem states that in this case each Hodge
manifold is projective algebraic with $N\omega-[H]=d\alpha$, where $H$ is a hyperplane section,
$N$ a positive integer, and $\alpha$  a current of degree 1. Note that $[H]$ is a $*\omega$-submanifold
not an $\omega$-submanifold.

\medskip\noindent
{\bf The Hodge Question for the Class of $\phi$}. 
Suppose $(X,\phi)$ is a Hodge manifold.  When does there exist a $*\phi$-cycle $T$ cohomologous
to $N\phi$ for some positive integer $N$, i.e.,
$$
N\phi-T\ =\ d\alpha
\eqno{(6.1)}
$$
for some current $\alpha$ of degree $p-1$?
\medskip

Recall that by definition a  $*\phi$-cycle is automatically $*\phi$-positive, so this is, more precisely,  the ``Hodge Question with Positivity'' for $\phi$.

\Remark{}  If   equation (6.1) (called the {\sl spark equation}) has a solution, then $\alpha$ determines a differential character on $X$.  (See [HLZ]  for more details.)

\Ex{6.3}  In [L$_3$] an  example is constructed of a Hodge manifold $(X,\phi)$ for which no such cycle exists.
More specifically, a parallel self-dual 4-form $\phi$ of comass 1 is constructed  on a flat torus $X$ of dimension 8 with the property that  
$[\phi]\in {\widetilde H}^4(X,\bbz)$, but there exist no ${\phi}$-cycles  whatsoever on $X$.

\Ex{6.4} Consider  the fundamental bi-invariant 3-form $\O$ on a compact simple Lie group $G$,
 normalized to be the generator of $H^3(G,\bbz)\cong \bbz$. By   rescaling the bi-invariant metric on $G$ we may assume that $\O$ has comass 1, i.e., that 
 $(G,\O)$ is a Hodge manifold.  H. Tasaki [T, Thm. 7],
  building on work of  Dao \v Cong Thi [Thi], showed that
 $*\O$-cycles exist and so there is a positive answer to the Hodge Question for $\O$.  Later
 R. Bryant [B] proved  that $\O$ is  itself cohomologous to a $*\O$-cycle, i.e., we can take $N=1$
 in this case.  
 He also showed that all $*\O$-cycles are  sums of  singular semi-analytic subvarieties congruent to irreducible semi-analytic  components of the cut-locus of the exponential map.
\medskip

Recall from Definition 2.1 that
$$
\Lambda(\phi) \  \equiv \  {\rm span}\{ G(\phi)\}\ \subset\ \Lambda_p TX.
$$
\Def{6.5}  A $p$-dimensional current $T$, representable by integration,  is said to be of 
{\bf type} $\Lambda(\phi)$ if 
$\overrightarrow T_x \in \L(\phi) \subset \L_pT_xX$ for $\| T \|$-a.a. $x$.  
This definition   extends to arbitrary currents $T$ of dimension $p$. 
If $T(\psi)=0$ for all smooth $p$-forms $\psi$ such that $\psi\bigr|_{\L(\phi)}=0$,
then $T$ is said to be of  {\bf type} $\Lambda (\phi)$.
\medskip

Each $\phi$-cycle $T$ is of type $\Lambda (\phi)$ since ${\overrightarrow T}_x \in G(\phi)$
$\|T\|$-almost everywhere.

\Def{6.6}  A {\bf  $\L(\phi)$-cycle} is a $d$-closed, $p$-dimensional locally rectifiable current of 
type $\L(\phi)$.
\medskip

There is a natural necessary condition for a class $c\in {\widetilde H}_p $ to be represented by
a $\Lambda (\phi)$-cycle. In fact we have the following more general statement. 

\Prop{6.7}  {\sl  If a class $c\in \IH_p(X,\bbz)$ is represented by a current of type $\L(\phi)$, then the harmonic representative of $c$ must be of type $\L(\phi)$.
}

\pf
Recall that the Hodge decomposition: 
$\ce^p(X) = \bbh^p(X) \oplus {\rm Image }(d)\oplus {\rm Image }(d^*)$, is a $C^\infty$-decomposition, and therefore induces a corresponding decomposition of currents:
${\ce'}_p(X) = \bbh_p(X) \oplus {\rm Image }(\partial)\oplus {\rm Image }(\partial^*)$.
Note that the orthogonal bundle projection $P_{\Lambda(\phi)} : \Lambda_pTX \to \Lambda(\phi)$ 
is a parallel operator. 
It was proved by Chern [Ch] that any such operator commutes with harmonic projection $\bbh$.
Suppose now that $c$ is represented by  a current $T$ of type $\L(\phi)$.  Then 
$P_{\Lambda(\phi)} (T) = T$ and therefore $P_{\Lambda(\phi)} \bbh(T)
= \bbh(T)$.\qed

\medskip\noindent
{\bf  I.  The Hodge Question.}  Suppose  $c   \in  \IH_p(X,\bbz)$ is a class whose harmonic representative
is of type $\L(\phi)$.  When does there exist an integer $N$ and a $\L(\phi)$-cycle $T$ with
$T\in Nc$?

\Remark{6.8} Example 6.3 above gives a parallel 
calibration $\phi$ on a flat 8-dimensional torus $X$ and an integral  class 
$c\in \widetilde {H}_4(X, \bbz)$
of type $\L(\phi)$ for which no such current exists. 

\Remark{6.9} The Hodge Question  is  a direct generalization of the
standard Hodge Conjecture for algebraic cycles on a complex projective manifold,
since we know from [HS],  [Sh] and [Alex]  that for 
$\phi = \omega^p/p!$ ($\omega$ =  the K\"ahler form), any 
$\L(\phi)$-cycle is an algebraic $p$-cycle.

\Remark{6.10} Any locally finite integer sum of $\phi$-cycles is a $\L(\phi)$-cycle.  
However, the converse is
completely open outside of the K\"ahler case.  Moreover, 
even though it holds in the K\"ahler case (cf. Remark 6.9),  there is  no proof of this fact by the standard methods of regularity
 in Geometric Measure Theory.

Before trying to prove that a general $\L(\phi)$-cycle is a sum of $\phi$-cycles, one
would like the calibration $\phi$ to have the following algebraic property (6.2).
 Equation (2.4) says that 
$
G(p, T_xX)\cap \L_+(\phi)\ =\ G(\phi)
$
so that $\phi$-cycles and $\lp$-cycles are the same thing.  Most parallel calibrations
(see [HL$_3$, p. 68])  are known to satisfy
$$
G(p, T_xX)\cap \L(\phi)\ =\ G(\phi) \cup (-G(\phi)).
\eqno{(6.2)}
$$
In this case $T$ is a $\lp$-cycle if and only if $\pm \overrightarrow {T_x} \in G_x(\phi)$ for
$\|T\|$-a.a.$\ x$.
Consequently, $T$ decomposes into $T^+-T^-$ with both $\overrightarrow {T_x}^{\pm}  \in G_x(\phi)$,
but, even in the K\"ahler case, one can not show directly that $T^+$ and $T^-$ are $d$-closed. 

\medskip

There are versions of the Hodge Question involving ``positivity'' which may have more hope.
For example:

\medskip\noindent
{\bf  II. The Hodge Question (with positivity).}
Suppose $c\in  \IH_p(X,\bbz)$ is a class whose harmonic representative is strictly $\lp$-positive.
When does there exist an integer $N$ and a $\phi$-cycle $T$ with $T\in Nc$?

\Remark{6.11}  If the current $*\phi$ (of dimension $p$)  is strictly $\lp$-positive, then for any form 
$\psi$ of type $\L(\phi)$, there exists an integer $\ell$ such that $\psi+\ell(*\phi)$ is strictly
$\lp$-positive.  This applies for example to the harmonic representative of $c$ in Hodge Question II.
Consequently, one can see that {\sl if $(X, *\phi)$ is a Hodge manifold with a solution to 
(6.2), then the Hodge Question I follows from II.}

\Remark{6.12}  The point of  Hodge Question II  is that one is asking for a $\phi$-cycle $T$.
These are   automatically $\lp$-positive and therefore satisfy the strong regularity Theorem 2.10.
\medskip

 Federer and Fleming [FF]  showed that each class $c\in H_p(X,\bbz)$ contains a rectifiable cycle
 $T$ with $M(T)\leq M(S)$ for all other rectifiable cycles $S\in c$.  
 Let   $\|c\|_\bbz =M(T)$ denote  this minimum.  
 Let   $\|\widetilde c\|_\bbr$ denote  the infimum
   of the masses $M(S)$ taken over {\sl all} closed currents
 homologous to $T$, i.e., over all currents in the real homology class $\wt c \in {\wt H}_p(X,\bbz)$
 corresponding to $c$.
 
 \Prop{6.13} {\sl  Fix $c\in H_p(X,\bbz)$ and suppose the corresponding class
 $\wt c \in  {\wt H}_p(X,\bbz)$ has a smooth representative $\psi$ which is 
 $\L_+(\phi)$-positive.  Then there exists a $\phi$-cycle $T\in c$ if and only if 
 $\|c\|_\bbz  = \|\wt c\|_\bbr$.}
  
 \pf    If $c$ contains a $\phi$-cycle $T$,    
 then  $\|c\|_\bbz \leq M(T) = T(\phi) = S(\phi) \leq M(S)$ for all currents 
 $S\in \wt c$. Hence, $\|c\|_\bbz  \leq \|\wt c\|_\bbr$, and the inequality $\|\wt c\|_\bbr  \leq \|c\|_\bbz$
 is clear from the definitions.
 
 Conversely, if $\|c\|_\bbz  = \|\wt c\|_\bbr$, then the Federer-Fleming solution $T$ (a rectifiable current)
 satisfies $M(T) =  \|c\|_\bbz  = \|\wt c\|_\bbr \leq M(\psi) =\int\psi\wedge\phi = T(\phi)$ 
 since $T$ and $\psi$ are homologous.  However, $T(\phi) \leq M(T)$ with 
 equality iff $T$ if a $\phi$-cycle. \qed
 
 \medskip
 
 Thus the positive Hodge question can be rephrased as follows:  When does there exist an integer
 $k$ with $\|kc\|_\bbz  = \|k\wt c\|_\bbr$?  
 
 \medskip
 \noindent
 {\bf Final Note:}  Since   $\|k\wt c\|_\bbr=k\|\wt c\|_\bbr$, we have $\|kc\|_\bbz  = \|k\wt c\|_\bbr$ iff
$ \|\wt c\|_\bbr =\smfrac 1 k \|kc\|_\bbz$.  Federer has shown in [F$_1$]  that
$$
 \|\wt c\|_\bbr \ =\ \lim_{k\to\infty} \smfrac 1 k  \|kc\|_\bbz
$$

\vfill \eject


\centerline{\bf  Appendix: The Reduced $\f$-Hessian.}\medskip

We assume throughout this section that $\L(\f)$ is a vector subbundle
of $\L_pTX$, and we let $\L(\f)\subset \L^p T^*X$ denote the corresponding bundle under
the metric equivalence $\L_pTX\cong \L^p T^*X$.

\Def{A.1}  The {\bf reduced $\f$-hessian} $\RH : C^\infty(X)\to 
\G(X, \L(\f))$ is defined to be $\ch^\f$ followed by orthogonal projection  onto the 
subbundle  $\L(\f)\subset \L^p T^*X$.  \medskip

Note that a function $f$ is $\f$-pluriharmonic if and only if $\RH (f)=0$.

Note also that if $\f$ is parallel, then $\RH = \overline d d^\f$ where $ \overline d $ denotes
the exterior derivative followed by orthogonal projection onto $\L(\f)$.

For most of the calibrations considered as examples in this paper, the image of the map
$\BM_\f:\Sym^2(TX)\to \L^p T^*X$ is contained in $\L(\f)$, or equivalently, $\RH = \ch^\f$. 
For reference, $\RH = \ch^\f$ in the following cases.
 \smallskip

(1) \ \ $\f = {1\over p!} \o^p$, the $p$th power of the K\"ahler form,

(2) \ \ $\f$ Special Lagrangian

(3) \ \ $\f$ Associative, Coassociative or Cayley

(4) \ \ $\f$ the fundamental 3-form on a simple Lie group.

\smallskip

Exceptions will be discussed at the end of this appendix.

Even when  $\RH = \ch^\f$ the following proposition is important.

\Prop{A.2}  {\sl  Suppose $f$ is a distribution on $X$.  Then $f$ is \fp if and only if 
$\RH(f)\equiv R$ is representable by integration and $\overrightarrow R\in \L^+(\f)$
$\|R\|$-a.e., that is, if and only if $\RH(f)$ is a $\L^+(\f)$-positive current.}

\medskip
The proof is similar to the proof of Theorem \DD.14 and is omitted.

\Def{A.3}  The $\f$-Grassmannian $G(\f)$ {\bf involves all the variables}  if, for $u\in TX$, 
the condition $u\hk \x=0$ for all $\x \in G(\f) $ implies   $u=0$
 
 \Ex{} The 2-form $\f\equiv dx_1\wedge dx_2+\lambda dx_3 \wedge dx_4$ with $|\lambda|<1$
 is a calibration on $\bbr^4$ which involves all the variables (See Section \AA), but the only
 $\x\in G(\f)$ is the $x_1,x_2$ plane so that $G(\f)$ does not involve all the variables.
 
 \Prop{A.4} {\sl  The operator $\RH$ is over-determined elliptic  if and only if $G(\f)$ involves all the variables.}
 
 \pf
 We  need only consider the case $\f\in \L^p V$, where $V$ is an inner product space.
 The symbol of $\ch^\f$  at $u\in V$ is $u\wedge(u\hk \f)$.  Hence, the reduced 
 operator $\RH$ is elliptic if and only if 
 $$
( u\wedge(u\hk \f))(\x)=0\ \  \forall \ \ \x \in G(\f) \qquad\Rightarrow \qquad u=0.
 $$
 For $\x\in G(p,V)$ and $u\in V$, let $u=a+b$ with $a\in\span \x$ and $b\perp\span \x$.
 Then $(u\wedge(u\hk \f))(\x) = \f(u\wedge(u\hk \x)) = \f((a+b)\wedge( a\hk \x))
 = |a|^2\f(\x) +\f(b\wedge(a\hk\x))$.  If $\x\in G(\f)$, then
 $\f(b\wedge(a\hk\x))=0$ by the First Cousin Principle [HL$_5$, Lemma 2.4] , and  $\f(\x)=1$.  Hence, 
 $(u\wedge(u\hk \f))(\x) = |a|^2=0$ if and only if $u\hk\x=0$. \qed
 
 \medskip
 
 One can easily reduce a calibration to the elliptic case.
 
 \Prop{A.5}  {\sl  Suppose $\f\in \L^p V$ is a calibration.  Define $W\subset V$ by 
 $$
 W^\perp\ \equiv\ \bigcap_{\x\in G(\f)}
(\span\x)^\perp 
$$
and set $\psi = \phi\bigr|_W$.  Then $\psi\in \L^pW$ is a calibration and $G(\psi)$ involves
all the variables in $W$.  Moreover, $G(\f)=G(\psi)$ and the reduced operators
$\RH$ and $\overline {\ch}^{\psi}$ agree.}
 
 \pf
 Obviously $\psi$ is a calibration and $G(\psi)\subset G(\f)$.  
 If $\x\in G(\f)$, then $\span \x\subset W$ and hence $\f(\x) = \psi(\x)$. Thus $G(\f)=G(\psi)$.
 By construction $G(\psi)$ involves all the variables in $W$.  Finally, for all $\x \in G(\f)$,
we have $\RH (f)(\x) =\tr_\x \Hess f= \overline {\ch}^{\psi}(f)(\x)$.\qed
 
 \Ex{} Let  $\Psi\in \L^4_{\bbr}\bbh^n$ be the quaternionic calibration on $\bbh^n$.
 One can show that $dd^\Psi f=0$ if and only if $\Hess f=0$.
 However,
  $$
\overline{d} d^\Psi f\ =\  \overline {\ch}^{\Psi}(f)\  =\ 
 \BM_{\Psi}\left(\!\!\left({\partial^2f\over \partial q_{\a} \partial \bar{q}_{\beta}}\right)\!\!\right),
 $$
 that is, the reduced hessian is isomorphic to the quaternionic hessian
 $
 \left({\partial^2f\over \partial q_{\a} \partial \bar{q}_{\beta}}\right)
 $

\vfill\eject



\centerline{\bf References}

\vskip .2in

\noindent
[Al]   S. Alesker,  {\sl  Non-commutative linear algebra and  plurisubharmonic functions  of quaternionic variables}, Bull.  Sci.  Math., {\bf 127} (2003), 1-35. also ArXiv:math.CV/0104209.  

\smallskip

\noindent
[AV]   S. Alesker and M. Verbitsky,  {\sl  Plurisubharmonic functions  on hypercomplex manifolds and HKT-geometry}, J. Geom. Anal. {\bf 16} (2006), 375-399.  arXiv: math.CV/0510140

\smallskip

\noindent
[Alex]   H. Alexander,  {\sl  Holomorphic chains and the support hypothesis conjecture}, J. Amer. Math.
Soc., {\bf 10} (1997), 123-138.

\smallskip

\noindent
[A] F. J. Almgren, Jr.,  {\sl  $Q$-valued functions minimizing Dirichlet's integral and the regularity of area-minimizing rectifiable currents up to codimension 2},  World Scientific Monograph Series
in Mathematics, 1, World Scientific Publishing Co.River Edge, NJ, 2000.

\smallskip

\noindent
[BJS]   L. Bers, F. John and M. Schechter,  {Partial Differential Equations}, Interscience, J. Wiley,
1964.

\smallskip

\noindent
[Be]   A. Besse,  {Einstein Manifolds}, Springer-Verlag, New York, 1987.
1987.

\smallskip

\noindent
 [B]   R. L.  Bryant,  {\sl
 Calibrated cycles of codimension 3 in compact simple Lie groups},
  to  appear.
 \smallskip

\noindent
 [BS]   R. L.  Bryant and S. M. Salamon,  {\sl
 On the construction of some complete metrics with exceptional holomony},
  Duke Math. J. {\bf 58} (1989),   829-850.

 \smallskip

\noindent
 [C]   E.  Calabi,  {\sl
 M\'etriques k\"ah\'eriennes et fibr\'es holomorphes},
 Annales scientifiques de l'\'Ecole Normale Superieure {\bf 12} (1979),   269-294.

 \smallskip

\noindent
 [Ch]   S.-S. Chern,  {\sl
 On a generalization of K\"ahler geometry},  Lefschetz Jubilee Volume, Princeton Univ. Press,
 (1957), 103-121.

 \smallskip

 \noindent
[DS] J. Duval and N. Sibony, {\sl
Polynomial convexity, rational convexity and currents},
  Duke Math. J. {\bf 79}  (1995),     487-513.

 \smallskip

\noindent
 [C]   T. Eguchi and A. J. Hanson,  {\sl
Asymptotically flat solutions to Euclidean gravity},
Physics Letters  {\bf 74B} (1978),   249-251.

 \smallskip

\noindent
[F$_1$]   H. Federer, Geometric Measure  Theory,
 Springer--Verlag, New York, 1969.

 \smallskip

\noindent
 [F$_2 $]   H. Federer, {\sl  Real flat chains, cochains and variational problems},
Indiana Univ.  Math. J.   {\bf 24}    (1974/75),   351-407.

 \smallskip

\noindent
 [FF]   H. Federer and W. Fleming, {\sl
Normal and Integral Currents},
Annals of Math.  {\bf 72} (1960),   458-520.

 \smallskip

\noindent
 [Fu]   L. Fu, {\sl  On the boundaries of Special Lagrangian submanifolds},
Duke Math. J.   {\bf 79}   no. 2 (1995),   405-422.

 \smallskip

\noindent
 [GL]  K. Galicki and B. Lawson, {\sl  Quaternionic reduction and quaternionic orbifolds}, Math. Ann. {\bf 282} (1989), 1-21.

 \smallskip

\noindent
[GR]   H. Grauer and R. Remmert,  Coherent Analytic Sheaves, Springer-Verlag, Berlin-Heidelberg, 1984.
\smallskip

\noindent
[GZ]  V. Guedj and A. Zeriahi,     
{\sl    Intrinsic capacities on compact K\"ahler manifolds},    
Preprint Univ. de Toulouse , 2003

\smallskip

\noindent
[H$_1$]  F.R. Harvey,
Holomorphic chains and their boundaries, pp. 309-382 in ``Several Complex
Variables, Proc. of Symposia in Pure Mathematics XXX Part 1'', 
A.M.S., Providence, RI, 1977.

\noindent
[H$_2$]  F.R. Harvey,
Spinors and Calibrations,  Perspectives in Mathematics, vol. 9 Academic Press, Boston, 1990

\noindent
[HK] F. R. Harvey and  A.  W. Knapp, {\sl  Positive (p,p)-forms, Wirtinger's inequality and currents}, 
Value-Distribution Theory, Part A (Proc. Tulane Univ. Program on Value-Distribution Theory
in Complex Analysis and Related Topics in Differential Geometry, 1972-73),  pp. 43-62,
Dekker, New York, 1974.

 \smallskip

\noindent
[HL$_1$] F. R. Harvey and H. B. Lawson, Jr, {\sl On boundaries of complex
analytic varieties, I}, Annals of Mathematics {\bf 102} (1975),  223-290.

 \smallskip

\noindent
[HL$_2$] F. R. Harvey and H. B. Lawson, Jr, {\sl On boundaries of complex
analytic varieties, II},  Annals of Mathematics {\bf 106} (1977), 
213-238.

 \smallskip

   \noindent 
 {[HL$_3$]} F. R. Harvey and H. B. Lawson, Jr, {\sl Calibrated geometries},  Acta Mathematica 
{\bf 148} (1982), 47-157.

 \smallskip

 
   \noindent 
 {[HL$_4$]} F. R. Harvey and H. B. Lawson, Jr, {\sl Boundaries of positive holomorphic chains
 and the relative Hodge question}, Stony Brook Prerprint, 2005. ArXiv:math/0610533

 \smallskip

   \noindent 
 {[HL$_5$]} F. R. Harvey and H. B. Lawson, Jr, {\sl An introduction to potential theory in calibrated geometry}, 
 Amer. J. Math. (to appear).  ArXiv:0710.3920

 \smallskip

   \noindent 
 {[HL$_6$]} F. R. Harvey and H. B. Lawson, Jr, {\sl  Plurisubharmonicity in a general geometric context}, Stony Brook Prerprint, 2008.

 \smallskip

   \noindent 
 {[HL$_7$]} F. R. Harvey and H. B. Lawson, Jr, {\sl  Dirichlet Duality and the Nonlinear Dirichlet Problem},
 Comm. Pure and Appl. Math.   {\bf 62} (2009), 396-443.  ArXiv:0710.3991

 \smallskip

\noindent
[HLZ] F. R. Harvey, H. B. Lawson, Jr. and J. Zweck, {\sl A
deRham-Federer theory of differential characters and character duality},
Amer. J. of Math.  {\bf 125} (2003), 791-847. ArXiv:math.DG/0512251

 \smallskip

\noindent
[HP] F. R. Harvey, J. Polking, {\sl Extending analytic objects},
Comm.  Pure Appl. Math. {\bf 28} (1975), 701-727.

 \smallskip

\noindent
[HW$_1$] F. R. Harvey,  R. O. Wells, Jr.,  {\sl Holomorphic approximation and hyperfunction 
theory on a $C^1$ totally real submanifold of a complex manifold},
  Math.  Ann. {\bf 197} (1972),  287-318.

 \smallskip

\noindent
[HW$_2$] F. R. Harvey,  R. O. Wells, Jr.,  {\sl Zero sets of non-negatively strictly plurisubharmonic
functions},
  Math.  Ann. {\bf 201} (1973),  165-170.

 \smallskip

   \noindent
[HS]    F.R. Harvey and B. Shiffman,    {\sl  A characterization of
holomorphic chains},    Ann. of Math.,
 {\bf 99}  (1974), 553-587.

\smallskip

   \noindent
[Ha$_1$]    M. Haskins,    {\sl  Special Lagrangian cones},    
Amer. J.  Math.,
 {\bf 126}  (2004), 845-871.

\smallskip

   \noindent
[Ha$_2$]    M. Haskins,    {\sl  The geometric complexity of special Lagrangian $T^2$-cones},    
Inventiones Math.,
 {\bf 157}  (2004), 11-70.

\smallskip

   \noindent
[HaK]    M. Haskins and N. Kapouleas,    {\sl  Special Lagrangian cones with higher genus links}.  ArXiv:math.DG/0512178.

\smallskip

   \noindent
[HaP]    M. Haskins and T. Pacini,    {\sl  Obstructions to special Lagrangian desingularizations,
and the Lagrangian prescribed boundary problem},    
to appear in Geometry and Topology.  ArXiv:math.DG/0609352.

\smallskip

\noindent
[J$_1$]   D. D. Joyce,     
{    Compact Manifolds with Special Holonomy},    
Oxford University Press, Oxford, 2000.

\smallskip

\noindent
[J$_2$]   D. D. Joyce,     
{    Special Lagrangian submanifolds with isolated conical siingularities, I. Regularity},    
 Ann. Global Anal. Geom. {\bf  25} (2004), 201-251. ArXiv:math.DG/0211294.

\smallskip

\noindent
[J$_3$]   D. D. Joyce,     
{    Special Lagrangian submanifolds with isolated conical siingularities, II. Moduli spaces},    
 Ann. Global Anal. Geom. {\bf  25} (2004), 301-352.  ArXiv:math.DG/0211295.

\smallskip

\noindent
[J$_4$]   D. D. Joyce,     
{    Special Lagrangian submanifolds with isolated conical siingularities, III. Desingularization,
the unobsructed case},    
to appear in Ann. Global Anal. Geom., ArXiv:math.DG/0302355.

\smallskip

\noindent
[J$_5$]   D. D. Joyce,     
{    Special Lagrangian submanifolds with isolated conical siingularities, IV. Desingularization,
 obsructions and families},    
to appear in Ann. Global Anal. Geom., ArXiv:math.DG/0302356.

\smallskip

\noindent
[J$_6$]   D. D. Joyce,     
{    Special Lagrangian submanifolds with isolated conical siingularities, V. Survey and applications},    
J. Diff. Geom.  {\bf 63} (2003), 279-347.  

\smallskip

\noindent
[J$_7$]   D. D. Joyce,     
{    U(1)-invariant special Lagrangian 3-folds, II. Existence of singular solutions},    
Adv. in Math.  {\bf 192} (2005), 72-134.  

\smallskip

   \noindent
[K]    J. King,    {\sl   The currents defined by analytic varieties},  Acta Math.   {\bf 127}  no. 3-4 (1971),
185-220.

\smallskip

   \noindent
[L$_1$]    H. Blaine Lawson, Jr.,      Minimal Varieties in Real and Complex Geometry, Les Presses de 
L'Universite de Montreal,  1974.

 \smallskip

   \noindent
[L$_2$]    H. Blaine Lawson, Jr.,    {\sl  Minimal Varieties},    
Proceedings of Symposia in Pure Mathematics 
{\bf 27} (1974), 61-93.

 \smallskip

   \noindent
[L$_3$]    H. Blaine Lawson, Jr.,    {\sl  The stable homology of a flat torus},    
Mathematica Scandinavica 
{\bf 36} (1975), 49-73.

\smallskip

 \noindent 
 {[L$_4$]}   H. Blaine Lawson, Jr.,    {\sl Geometric aspects of the genealized Plateau problem}, 
pp. 7-13 in  Proceedings  of the International Congress of Mathematicians, 1974,
vol.2, Vancouver, Canada, 1975.

 \smallskip

 \noindent 
 {[L$_5$]}   H. Blaine Lawson, Jr.,    {\sl The question of holomorphic carriers}, Proceedings of
Symposia in  Pure Mathematics {\bf 30}, American Mathematical 
Society (1976), 115-124.

 \smallskip

 \noindent 
 {[M]} J. Milnor, {Morse Theory}, Annals of Math. Studies no. {\bf 51}, Princeton University Press,
 Princeton, N.J.,  1963.
 \smallskip

   \noindent
[S]   H. H. Schaefer,  Topological Vector Spaces,    Springer Verlag,
New York,  1999.

\smallskip

   \noindent
[Sh]      B. Shiffman,    {\sl Complete characterization of holomorphic chains of codimension one},    
Math. Ann. {\bf 274} (1986), 233-256.
\smallskip

\noindent
[T] H. Tasaki,  {\sl  Certain minimal or homologically volume minimizing submanifolds in compact symmetric spaces},  Tsukuba J.  Math.  , {\bf 9} (1985), 117-131.

\smallskip

\noindent
[Thi]  Dao \v Cong Thi,  {\sl  Real minimal currents in compact Lie groups},  Trudy Sem. Vektor. Tenzor. Anal. , {No. 19} (1979), 112-129.

\smallskip

\noindent
[U]  I.  Unal, Ph.D. Thesis, Stony Brook, 2006.

   \noindent
[V]    M. Verbitsky.,    {\sl  Manifolds with parallel differential forms and K\"ahler identities
for $G_2$-manifolds},  arXiv : math.DG/0502540 (2005).
 \smallskip

\end